\documentclass[reqno,10pt]{article}
\usepackage{amssymb, amsmath}
\newcommand{\redch}[2]{{{{#2}}}}
\newtheorem{theorem}{Theorem}
\newtheorem{lemma}[theorem]{Lemma}
\newtheorem{proposition}[theorem]{Proposition}
\newtheorem{corollary}[theorem]{Corollary}

\newtheorem{remark}[theorem]{Remark}
\newenvironment{pf}{\medskip\par\noindent{\it Proof\ }:}{\hfill$\Box$\medskip\par}
\newenvironment{pf*}[1]{\medskip\par\noindent{\it Proof of #1\ }:}{\hfill$\Box$\medskip\par}
\def\eqref#1{$(\ref{#1})$}

\def\SET#1#2{\left\{\,#1\,:\,#2\,\right\}}
\def\SSET#1#2{\{\,#1\,:\,#2\,\}} 

\def\R{\mathbb{R}}
\def\N{\mathbb{N}}
\def\S{\mathbb{S}}

\def\MIN{\mbox{ \rm minimize }}

\def\ST{\mbox{ \rm subject to }}

\def\FEAS(#1){\mathop{\rm Feas}\langle #1 \rangle} 
\def\PROBLEM#1{\ifmmode{\langle #1 \rangle}\else\ifinner{\langle #1 \rangle}
        \else{$\langle #1 \rangle$}\fi\fi}
\def\mathbfit#1{\mbox{\mathversion{bold}$#1$\mathversion{normal}}}
\def\b0{\mathbfit{0}}

\begin{document}
\begin{center}
\sc\Large
A Perturbed Sums of Squares Theorem for Polynomial Optimization and its Applications
\end{center}
\bigskip
\begin{center}
Masakazu Muramatsu\footnote{Department of Communication Engineering and Informatics, The University of Electro-Communications, 
		1-5-1 Chofugaoka, Chofu-shi, Tokyo, 182-8585 JAPAN. 
		muramatu@cs.uec.ac.jp
}, Hayato Waki\footnote{Institute of Mathematics for Industry, Kyushu University, 744 Motooka, Nishi-ku, Fukuoka 819-0395, JAPAN. 
waki@imi.kyushu-u.ac.jp
}
and Levent Tun{\c c}el\footnote{
Department of Combinatorics and Optimization, 
Faculty of Mathematics, 
University of Waterloo, 
Waterloo, Ontario N2L 3G1
CANADA. 
ltuncel@math.uwaterloo.ca 
}
\end{center}
\begin{center}
\today
\end{center}
\begin{center}
\textbf{Abstract} \\
\smallskip
\begin{minipage}[t]{12cm}
We consider a property of positive polynomials on a compact set 
with a small perturbation.
When applied to a Polynomial Optimization Problem (POP), the property implies that 
the optimal value of the corresponding SemiDefinite Programming (SDP) relaxation with sufficiently 
large relaxation order is bounded from below by $(f^\ast - \epsilon)$ and 
from above by $f^\ast + \epsilon (n+1)$,
where $f^\ast$ is the optimal value of the POP. 
We propose new  SDP relaxations for POP based on modifications of existing sums-of-squares representation theorems. An advantage of our  SDP relaxations is that in many cases they are of considerably smaller dimension than those originally proposed by Lasserre. We present some applications and the results of our computational experiments. 
\end{minipage}
\end{center}
\noindent

\renewcommand{\theenumi}{\roman{enumi}}
\section{Introduction}\label{sec:introduction}
\subsection{Lasserre's SDP relaxation for POP}
We consider the POP:
\begin{equation} \label{eq:POP0}
 \MIN\ f(x) \ \ST\ f_i(x) \geq 0 \ (i=1,\ldots,m),
\end{equation}
where $f$, $f_1,\ldots,f_m: \R^n\rightarrow \R$ 
are polynomials.
The feasible region is denoted by 
$K = \SSET{x\in\R^n}{f_j(x)\geq 0  \ (j=1,\ldots,m)}$.
Then it is easy to see that the optimal value $f^\ast$ can be represented as
\[
 f^\ast = \sup\SET{\rho}{f(x)-\rho\geq 0 \ (\forall x\in K)}.
\]

First, we briefly describe the framework of the SDP relaxation method 
for POP \eqref{eq:POP0} proposed by Lasserre \cite{LAS01}.
See also \cite{PARRILO03}. We denote the set of polynomials and sums of squares by $\R[x]$ and 
$\Sigma$, respectively.
$\R[x]_{r}$ is the set of polynomials whose degree is less than 
or equal to $r$. 
We let $\Sigma_r = \Sigma \cap \R[x]_{2r}$.
We define the quadratic module generated by $f_1,\ldots,f_m$ as
\[
M(f_1,\ldots,f_m) = \SET{\sigma_0+\sum_{j=1}^m \sigma_j f_j}%
    {\sigma_0,\ldots,\sigma_m \in \Sigma}.
\]
The truncated quadratic module whose degree is less than or equal to 
$2r$ is defined by 
\[
M_r(f_1,\ldots,f_m) = \SET{\sigma_0+\sum_{i=1}^m \sigma_j f_j}%
    {\sigma_0\in \Sigma_r,\sigma_j\in \Sigma_{r_j}(j=1,\ldots,m)},
\]
where $r_j = r - \lceil\deg f_j / 2 \rceil$ for $j=1,\ldots,m$.

Replacing the condition that $f(x)-\rho$ is nonnegative 
by a relaxed condition that the polynomial is contained in 
$M_r(f_1,\ldots,f_m)$, we obtain the following SOS relaxation:
\begin{equation} \label{eq:SOS1}
 \rho_r = \sup\SET{\rho}{f(x) - \rho \in M_r(f_1,\ldots,f_m)}.
\end{equation}
Lasserre\cite{LAS01} showed that 
$\rho_r\rightarrow f^\ast$ as $r\rightarrow\infty$ 
if $M(f_1,\ldots,f_m)$ is \redch{A}{a}rchimedean.
See \cite{Laurent09, prestel} for a definition of \redch{A}{a}rchimedean.
{An easy way to ensure that $M(f_1,\ldots,f_m)$ is \redch{A}{a}rchimedean is to make sure that $M(f_1,\ldots,f_m)$ contains a representation of a ball of finite (but possibly very large) radius. }
In particular, we point out that 
when $M(f_1,\ldots,f_m)$ is \redch{A}{a}rchimedean, $K$ is compact.

The problem \eqref{eq:SOS1} can be encoded as an SDP problem.
Note that we can express a sum of squares $\sigma\in\Sigma_r$
by using a positive semidefinite matrix $X\in \S^{s(r)}_+$
as $\sigma(x) = u_r(x)^T X u_r(x)$,
where $s(r) = {{n+r}\choose {n}}$
and $u_r(x)$ is the monomial vector which contains all the 
monomials in $n$ variables up to and including degree $r$ with an appropriate order.
By using this relation,
the containment by $M_r(f_1, \ldots, f_m)$ constraints in \eqref{eq:SOS1}, {\it i.e.}, 
\[
 f - \rho = \sigma_0 + \sum_{j=1}^m \sigma_{j} f_j,
\]
can be transformed to linear equations 
involving  semidefinite matrix variables corresponding to 
$\sigma_0$ and $\sigma_j$'s.

Note that, in this paper, we neither assume that
$K$ is compact nor that $M(f_1,\ldots,f_m)$ is \redch{Archmedean}{archimedean}. 
Still, the framework of Lasserre's SDP relaxation described above can be 
applied to \eqref{eq:POP0}, although the good theoretical convergence
property may be lost.

\subsection{Problems in the SDP relaxation for POP}

Since POP is NP-hard, solving POP in practice 
is sometimes extremely difficult.
The SDP relaxation method described above also has some difficulty.
A major difficulty arises from the size of the SDP relaxation 
problem \eqref{eq:SOS1}. In fact, \eqref{eq:SOS1} contains 
${{n+2r}\choose {n}}$ variables and 
$s(r)\times s(r)$ matrix. 
When $n$ and/or $r$ get larger, solving \eqref{eq:SOS1}
can become just impossible.

To overcome this difficulty, several techniques, using 
sparsity of polynomials, are proposed. See, e.g., 
\cite{KOJ03a,Lasserre06b,Laurent09,Nie07,WAKI06}.
Based on the fact that most of the practical POPs are 
sparse in some sense, these techniques exploit special sparsity structure
of POPs to reduce the number of variables and the size of the matrix variable in  the SDP
\eqref{eq:SOS1}.
\redch{}{Recent work in this direction, e.g., 
\cite{GHASEMI1,GHASEMI2} also exploit special structure of POPs
to solve larger sized problems. Nie and Wang \cite{NIE} proposes a use of regularization method for solving SDP relaxation problems instead of primal-dual interior-point methods. }

Another problem with the SDP relaxation is that \eqref{eq:SOS1} 
is often ill-posed.
In \cite{Henrion05, WAKI09, WAKI11}, strange behaviors of SDP solvers are 
reported. Among them is that an SDP solver returns an `optimal' value 
of \eqref{eq:SOS1} which is significantly different from the true
optimal value without reporting any numerical errors.
Even more strange is that the returned value by the SDP solver is 
nothing but the real optimal value of the POP \eqref{eq:POP0}.
We refer to this as  a `super-accurate' property of the SDP relaxation for POP.

\subsection{Contribution of this paper}\label{subsec:contribution}
POP contains very hard problems as well as some easier ones. We would like an approach which will exploit the structure in the easier instances of POP. In the context of current paper the notion of ``easiness" will be based on sums of squares certificate and sparsity. 
Based on Theorems \ref{th:main}, \ref{th:conv} and its variants, we propose new SDP relaxations. We call it \textit{Adaptive SOS relaxation} in this paper.
Adaptive SOS relaxations can be interpreted as relaxations of those originally proposed by Lasserre. As a result, the bounds generated by our approach cannot be superior to those generated by Lasserre's approach for the same order 
relaxations. However, Adaptive SOS relaxations are {of} significantly smaller \redch{dimension}{dimensions} (compared to Lasserre's SDP relaxations) and as the computational experiments in Section \ref{sec:experiment} indicate, we obtain very significant speed-up factors and we are able to solve larger instances and higher-order SDP relaxations. Moreover, in most cases, the amount of loss in the quality of bounds is small, even for the same order SDP relaxations.

\medskip

The rest of this paper is organized as follows.
Section \ref{sec:sdprelax} gives our main results and Adaptive SOS relaxation based on Theorem \ref{th:main}. 
In Section \ref{sec:experiment}, we present the results of some numerical experiments. We give a proof of Theorem \ref{th:main} and some of extensions, and the related work to Theorem \ref{th:main}  in Section \ref{sec:extensions}. 

\section{Adaptive SOS relaxation}\label{sec:sdprelax}

\subsection{Main results}\label{subsec:main}
We assume that there exists an optimal solution $x^\ast$ of \eqref{eq:POP0}.
Let 
\begin{eqnarray*}
 b &=& \max\left(1, \max\SSET{|x^\ast_i|}{i=1,\ldots,n}\right) \\
 B &=&[-b, b]^n.
\end{eqnarray*}
Obviously $x^\ast \in B$.
We define:
\begin{eqnarray*}
 \bar{K} &=& B\cap K \\
 R_j &=& \max\SET{|f_j(x)|}{x\in B} \ (j=1,\ldots,m). 
\end{eqnarray*}
Define also, for a positive integer $r$, 
\begin{eqnarray*}
 \psi_r(x) &=& - \sum_{j=1}^m f_j(x) \left(1-\frac{f_j(x)}{R_j}\right)^{2r}, \\
  \Theta_r(x) &=& 1 + \sum_{i=1}^n x_i^{2r}, \\
 \Theta_{r, b}(x) &=& 1 + \sum_{i=1}^n \left(\frac{x_i}{b}\right)^{2r}.
\end{eqnarray*}

We start with the following theorem.
\begin{theorem} \label{th:main}
Suppose that for $\rho\in \R$, $f(x) -\rho > 0$ for every $x\in \bar{K}$,
{\it i.e.}, $\rho$ is a lower bound of $f^\ast$.
\begin{enumerate}
\item Then there exists $\tilde{r}\in\mathbb{N}$ such that for all $r\ge\tilde{r}$, $f-\rho +  \psi_{r}$ is positive over $B$. 
\item In addition, for every $\epsilon>0$, 
there exists a positive integer $\hat{r}$ such that,
for every $r\geq\hat{r}$, 
\[
 f -\rho  + \epsilon \Theta_{r, b} + \psi_{\tilde{r}} \in \Sigma.
\]
\end{enumerate}
\end{theorem}

Theorem \ref{th:main} will be proved in Section \ref{sec:extensions}
as a corollary of Theorem \ref{th:sparsemain}. 
We remark that $\hat{r}$ depends on $\rho$ and $\epsilon$, while $\tilde{r}$ depends on $\rho$, but not $\epsilon$. The implication of this theorem is twofold. 

First, it elucidates the super-accurate property of 
the SDP relaxation for POPs.
Notice that by construction, 
$-\psi_{\tilde{r}}(x) \in M_{\bar{r}}(f_1,\ldots,f_m)$
where $\bar{r} = \tilde{r}\max_j(\deg(f_j))$.
Now assume that in \eqref{eq:SOS1}, $r\geq\bar{r}$.
Then, for any lower bound $\bar\rho$ of $f^\ast$,
Theorem \ref{th:main} means that
$f -\bar\rho + \epsilon\Theta_{r, b}\in M_{r}(f_1,\ldots,f_m)$
for arbitrarily small $\epsilon>0$ and sufficiently large $r$.

\redch{}{
Let us discuss this in more details.
Define $\Pi$ be the set of the polynomials such that
abosolute value of each coefficient is less than or equal to $1$.
Suppose that $\bar\rho$ is a ``close'' lower bound of $f^\ast$ such that
the system
$
 f - \bar{\rho} + \psi_{\tilde{r}} \in \Sigma
$ is infeasible.
Let us \textit{admit} an error $\epsilon$ in the above system, i.e., 
consider
\begin{equation} \label{eq:h}
 f - \bar{\rho} + \epsilon h + \psi_{\tilde{r}} \in \Sigma,
 \ h \in \Pi.
\end{equation}
The system \eqref{eq:h} restricts
the amount of the infinity norm error 
in the equality condition of the SDP relaxation problem
to be less than or equal to $\epsilon$.
Since we can decompose 
$h = h_+ - h_-$ where $h_+, h_-\in\Sigma\cap\Pi$,
now the system \eqref{eq:h} is equivalent with:
\begin{equation} \label{eq:hplus}
 f - \bar{\rho} + \epsilon h_+ + \psi_{\tilde{r}} \in \Sigma,
 \ h_+ \in \Pi\cap\Sigma.
\end{equation}
This observation shows that $-h_-$ is not the direction of errors.
Furthermore, because $\Theta_{r,b}\in\Pi\cap\Sigma$, the system 
\eqref{eq:hplus} is feasible due to ii of Theorem \ref{th:main}.
Therefore, if we admit an error $\epsilon$, 
the system $f- \bar\rho + \psi_{\tilde{r}} \in \Sigma$ is 
considered to be feasible,
and $\bar\rho$ is recognized as a lower bound for $f^\ast$.
}
\redch{}{As a result, we may obtain $f^\ast$ due to the 
numerical errors.}

\redch{}{
On the other hand, we point out that
when we do not admit an error, but are given 
a direction of error $h$ implicitly by the floating point arithmetic,
it does not necessarily satisfy the left inclusion of \eqref{eq:h}.
However, some numerical experiments show that this is true in most cases
(e.g., \cite{WAKI09}). The reason is not clear. 
}


Second, we can use the result to construct 
new sparse SDP relaxations for POP \eqref{eq:POP0}. Our SDP relaxation is weaker than Lasserre's, but the size of our SDP relaxation can become smaller than Lasserre's. As a result, for some large-scale and middle-scale POPs, our SDP relaxation can often obtain a lower bound, while Lasserre's cannot. 

A naive idea is that we use \eqref{eq:POP0} as is. 
Note that $-\psi_{\tilde{r}}(x)$
contains only monomials whose exponents are contained in 
\[
\bigcup_{j=1}^m \left(\mathcal{F}_j + \underbrace{\tilde{\mathcal{F}}_j+\cdots+\tilde{\mathcal{F}}_j}_{2\tilde{r}}\right), 
\]
where $\mathcal{F}_j$ is the \textit{support} of the polynomial $f_j$, {\it i.e.}, the set of exponents of monomials with nonzero coefficients in $f_j$, and $\tilde{\mathcal{F}}_j = \mathcal{F}_j\cup\{0\}$. 
To state the idea more precisely, we introduce some notation.
For a finite set $\mathcal{F}\subseteq \N^n$ and a positive integer $r$,
we denote $r\mathcal{F} = \underbrace{\mathcal{F}+\cdots+\mathcal{F}}_{r}$ and 
\[
 \Sigma(\mathcal{F}) = \SET{\sum_{k=1}^q g_k(x)^2}{\mbox{supp}(g_k) \subseteq \mathcal{F}}, 
\]
where $\mbox{supp}(g_k)$ is the support of $g_k$. 
Note that $\Sigma(\mathcal{F})$ is the set of sums of squares of polynomials
whose supports are contained in $\mathcal{F}$.

Now, fix an admissible error $\epsilon>0$ and $\tilde{r}$ as in Theorem \ref{th:main},  and consider: 
\begin{equation} \label{eq:sparse1}
\hat\rho(\epsilon, \tilde{r}, r) = 
 \sup\SET{\rho}{f-\rho + \epsilon \Theta_{r, b}- \sum_{j=1}^m f_j\sigma_j =
 \sigma_0, 
 \sigma_0\in \Sigma_r, 
 \sigma_j\in \Sigma(\tilde{r}\tilde{\mathcal{F}}_j)} 
\end{equation}
for some $r\geq\tilde{r}$.
Due to Theorem \ref{th:main}, \eqref{eq:sparse1}
has a feasible solution for all sufficiently large  $r$. 
\begin{theorem}\label{th:conv}
For every $\epsilon>0$, there exist $\tilde{r}, r\in\mathbb{N}$
such that 
$f^*-\epsilon\le \hat{\rho}(\epsilon, \tilde{r}, r)\le f^* + \epsilon(n+1)$. 
\end{theorem}
\begin{pf}
We apply Theorem \ref{th:main} to POP (\ref{eq:POP0}) with $\rho = f^*-\epsilon$. Then for any $\epsilon >0$, there exist $\hat{r}, \tilde{r}\in\mathbb{N}$ such that for every $r\ge \hat{r}$,  $f -(f^*-\epsilon) +\epsilon\Theta_{r, b} +\psi_{\tilde{r}}\in\Sigma$. 
Choose a positive integer $r\ge\hat{r}$ which satisfies 
\begin{equation}\label{choiceR}
r\ge\max\{
\lceil\deg(f)/2\rceil, \lceil(\tilde{r}+1/2)\deg(f_1)\rceil, \ldots, \lceil(\tilde{r}+1/2)\deg(f_m)\rceil
\}. 
\end{equation}
Then there exists $\tilde{\sigma}_0\in\Sigma_r$ such that $f -(f^*-\epsilon) +\epsilon\Theta_{r, b} +\psi_{\tilde{r}} = \tilde{\sigma}_0$,  because the degree of the polynomial in the left hand side is equal to $2r$. We denote $\tilde{\sigma}_j :=\left(1-f_j/R_j\right)^{2\tilde{r}}$ for all $j$. The triplet $(f^*-\epsilon, \tilde{\sigma}_0, \tilde{\sigma}_j)$ is feasible in (\ref{eq:sparse1}) because $\left(1-f_j/R_j\right)^{2\tilde{r}}\in\Sigma(\tilde{r}\tilde{\mathcal{F}}_j)$. Therefore, we have $f^*-\epsilon\le \hat{\rho}(\epsilon, \tilde{r}, r)$.   

We prove that $\hat{\rho}(\epsilon, \tilde{r}, r)\le f^* +\epsilon(n+1)$. We choose $r$ as in (\ref{choiceR}) and  consider the following POP:
\begin{equation}\label{origPOP}
\tilde{f}:=\inf_{x\in\mathbb{R}^n}\left\{
f(x) +\epsilon\Theta_{r, b}(x): 
f_1(x)\ge 0, \ldots, f_m(x)\ge 0
\right\}. 
\end{equation}
Applying Lasserre's SDP relaxation with relaxation order $r$ to (\ref{origPOP}), we obtain the following SOS relaxation problem:
\begin{equation} \label{tempSDPrelax}
\hat{\rho}(\epsilon, r):= \sup\SET{\rho}{f-\rho + \epsilon \Theta_{r, b} =
 \sigma_0 + \sum_{j=1}^m f_j\sigma_j, 
 \sigma_0\in \Sigma_{r}, 
 \sigma_j\in \Sigma_{r_j}}, 
\end{equation}
where $r_j := r -\lceil\deg(f_j)/2\rceil$ for $j=1, \ldots, m$. 
Then we have $\hat{\rho}(\epsilon, r)\ge \hat{\rho}(\epsilon, \tilde{r}, r)$ because $\Sigma(\tilde{r}\tilde{\mathcal{F}}_j)\subseteq\Sigma_{r_j}$ for all $j$. Indeed, it follows from (\ref{choiceR}) and the definition of $r_j$ that $r_j\ge \tilde{r}\deg(f_j)$, and thus $\Sigma(\tilde{r}\tilde{\mathcal{F}}_j)\subseteq\Sigma_{r_j}$.  

Every optimal solution $x^*$ of POP (\ref{eq:POP0}) is feasible for (\ref{origPOP}) and its  objective value is $f^* + \Theta_{r, b}(x^*)$. We have $f^* + \Theta_{r, b}(x^*)\ge \hat{\rho}(\epsilon, r)$ because (\ref{tempSDPrelax}) is the relaxation problem of (\ref{origPOP}). In addition, it follows from $x^*\in B$ that $n+1\ge \Theta_{r, b}(x^*)$, and thus $\hat{\rho}(\epsilon, \tilde{r}, r)\le \hat{\rho}(\epsilon, r)\le f^* + \epsilon(n+1)$. 
\end{pf}

Lasserre \cite{LAS01} proved the convergence of his  SDP relaxation under the assumption that the quadratic module $M(f_1, \ldots, f_m)$ associated with POP (\ref{eq:POP0}) is \redch{A}{a}rchimedean. In contrast, Theorem \ref{th:conv} does not require such an assumption and ensures that  we can obtain a sufficiently close {approximation} to the optimal value $f^*$ of POP (\ref{eq:POP0}) by solving (\ref{eq:sparse1}). 

We delete the perturbed part $\epsilon\Theta_{r, b}(x)$ from the above sparse relaxation  \eqref{eq:sparse1} in our computations, because it may be implicitly introduced in the computation by using floating-point arithmetic. In the above sparse relaxation  \eqref{eq:sparse1}, we have to consider only those positive semidefinite matrices 
whose rows and columns correspond to $\tilde{r}\tilde{\mathcal{F}}_j$ for $f_j$.
In contrast, in Lasserre's SDP relaxation, 
we have to consider the whole 
set of monomials whose degree is less than or equal to $r_j$
for each polynomial $f_j$.
Only $\sigma_0$ is large; it contains the set of all monomials
whose degree is less than or equal to $r$. 
However, since the other polynomials do not contain 
most of the monomials of $\sigma_0$, such monomials can  
safely be eliminated to reduce the size of $\sigma_0$ (as in \cite{KOJ03a}).
As a result, our sparse relaxation 
reduces the size of the matrix significantly if 
each $|\mathcal{F}_j|$ is small enough.
We note that in many of the practical cases, this in fact  is true.
We will call this new relaxation {\it Adaptive SOS relaxation} in the following.

\subsection{Proposed approach: Adaptive SOS relaxation}\label{subsec:overview}

\redch{A}{An} SOS  relaxation (\ref{eq:sparse1}) for POP (\ref{eq:POP0})  has been introduced. However, this relaxation has some weak points.
In particular, we do not know the value $\tilde{r}$ in advance.
Also, introducing small perturbation $\epsilon$ intentionally
may lead numerical difficulty in solving SDP.

To overcome these difficulties, we 
ignore the perturbation part $\epsilon\Theta_{r, b}(x)$ 
in (\ref{eq:sparse1}) because the perturbation part may be implicitly introduced by floating point arithmetic. 
In addition, we choose a positive integer $r$ and 
find $\tilde{r}$ by increasing $r$. 
Furthermore, we replace $\sigma_j\in\Sigma(\tilde{r}\tilde{\mathcal{F}}_j)$ by $\sigma_j\in\Sigma(\tilde{r}_j\tilde{\mathcal{F}}_j)$ in (\ref{eq:sparse1}), where $\tilde{r}_j$ is defined for a given integer $r$ as
\[
\tilde{r}_j = \left\lfloor
\frac{r}{\deg(f_j)} - \frac{1}{2}
\right\rfloor,
\]
to have $\deg(f_j\sigma_j) \le 2r$ for all $j=1, \ldots, m$. 
Then, we obtain the following SOS problem:
\begin{equation}\label{sosProb}
\rho^*(r) := \sup_{\rho\in\mathbb{R}, \sigma_0\in\Sigma_r, \sigma_j\in\Sigma(\tilde{r}_j\tilde{\mathcal{F}}_j)}\left\{
\rho: 
f - \rho -\sum_{j=1}^m f_j\sigma_j = \sigma_0
\right\}. 
\end{equation}
We call (\ref{sosProb}) {\it Adaptive SOS relaxation} for POP (\ref{eq:POP0}).  Note that we try to use numerical errors in a positive way;
even though Adaptive SOS relaxation has a different optimal 
value from that of POP, we may hope that the contaminated computation 
produces the correct optimal value of POP.

In general, we have $\Sigma(\tilde{r}_j\tilde{\mathcal{F}}_j)\subseteq \Sigma_{r_j}$ because of $\tilde{r}_j\deg(f_j)\le r_j$. Recall that $r_j=r -\lceil \deg(f_j)/2\rceil$ and is used in Lasserre's SDP relaxation (\ref{eq:SOS1}). This implies that Adaptive SOS relaxation is no stronger than  Lasserre's SDP relaxation, {\it i.e.}, the optimal value $\rho^*(r)$ is lower than or equal to the optimal value $\rho(r)$ of Lasserre's SDP relaxation for POP (\ref{eq:POP0}) for all $r$. We further remark that $\rho^*(r)$ may not converge \redch{}{to} the optimal value $f^*$ of POP (\ref{eq:POP0}). However, we can hope for the convergence of $\rho^*(r)$ to $f^*$ from Theorem \ref{th:main} and some numerical results in \cite{Henrion05, WAKI09, WAKI11}. 

\medskip

In the rest of this subsection, we provide a property of Adaptive SOS relaxation for the quadratic optimization problem
\begin{equation}\label{qop}
\inf_{x\in\mathbb{R}^n}\left\{
f(x) :=x^TP_0x + c_0^Tx :
f_j(x) := x^TP_jx + c_j^Tx + \gamma_j \ge 0 \ (j=1, \ldots, m)
\right\}. 
\end{equation}

The proposition implies that we do not need to compute $\rho^*(r)$ for even $r$.
\begin{proposition}
Assume that the degree $\deg(f_j) = 2$ for all $j=1, \ldots, m$ for QOP (\ref{qop}). Then, the optimal value $\rho^*(r)$ of Adaptive SOS relaxation is equal to $\rho^*(r-1)$ if $r$ is even.  
\end{proposition}
\begin{pf}
It follows from definition of $\tilde{r}_j$ that we have 
\[
\tilde{r}_j = \left\lfloor\frac{r-1}{2}\right\rfloor = \left\{
\begin{array}{cl}
\frac{r-1}{2} & \mbox{if }r \mbox{ is odd}, \\
 \frac{r}{2}-1 & \mbox{if }r \mbox{ is even}. 
\end{array}
\right.
\]
We assume that $r$ is even and give Adaptive SOS relaxation problems with relaxation order $r$ and $r-1$:
\begin{eqnarray}
\label{rho(r)}
\rho^*(r) &=& \sup\left\{
\rho :
\begin{array}{l}
f - \rho -\displaystyle\sum_{j=1}^m f_j\sigma_j = \sigma_0, \rho\in\mathbb{R}, \sigma_0\in\Sigma_r, \\
\sigma_j\in\Sigma\left(\displaystyle\left(\frac{r}{2}-1\right)\tilde{\mathcal{F}}_j\right)
\end{array}
\right\},\\
\label{rho(r-1)} 
\rho^*(r-1) &=& \sup\left\{
\rho :
\begin{array}{l}
f - \rho -\displaystyle\sum_{j=1}^m f_j\sigma_j = \sigma_0, \rho\in\mathbb{R}, \sigma_0\in\Sigma_{r-1}, \\
\sigma_j\in\Sigma\left(\left(\displaystyle\frac{r}{2}-1\right)\tilde{\mathcal{F}}_j\right)
\end{array}
\right\}. 
\end{eqnarray} 
We have $\rho^*(r) \ge \rho^*(r-1)$ for (\ref{rho(r)}) and (\ref{rho(r-1)}). All feasible solutions $(\rho, \sigma_0, \sigma_j)$ of (\ref{rho(r)}) satisfy the following identity:
\[
f_0 -\rho = \sigma_0 + \sum_{j=1}^m \sigma_j f_j. 
\]
Since $r$ is even, the degrees of $\sum_{j=1}^m \sigma_j(x) f_j(x)$  and $f_0(x)-\rho$ are  less than or equal to $2r-2$ and 2 respectively, and thus, the degree of $\sigma_0$ is less than or equal to $2r-2$. Indeed, we can write $\sigma_0(x) = \sum_{k=1}^{\ell}\left(g_k(x) + h_k(x)\right)^2$, where $\deg(g_k)\le r-1$ and $h_k$ is a homogenous polynomial with degree $r$. Then we obtain $0 = \sum_{k=1}^{\ell} h_k^2(x)$, which implies $h_k = 0$ for all $k=1, \ldots, \ell$. Therefore, all feasible solutions $(\rho, \sigma_0, \sigma_j)$ in SDP relaxation problem (\ref{rho(r)}) are also feasible in SDP relaxation problem (\ref{rho(r-1)}), and we have $\rho^*(r) = \rho^*(r-1)$ if $r$ is even.  
\end{pf}

\section{Numerical Experiments}\label{sec:experiment}

In this section, we compare 
Adaptive SOS relaxation with Lasserre's SDP relaxation and 
the sparse SDP relaxation using correlative sparsity proposed in \cite{WAKI06}. 
To this end, we  perform some numerical experiments. 
We observe from the results of our computational experiments that 
(i) 
although Adaptive SOS relaxation is often strictly weaker than Lasserre's, 
{\it i.e.}, the value obtained  by Adaptive SOS relaxation is less than Lasserre's, 
the difference is small in many cases, (ii) Adaptive SOS relaxation solves at least
10 times faster than Lasserre's in middle to large scale problems. 
Therefore, we conclude that Adaptive SOS relaxation can be more effective than Lasserre's for large- and middle-scale POPs. 
We will also observe a similar relationship against the sparse relaxation in \cite{WAKI06};
Adaptive SOS relaxation is weaker but much faster than the sparse one. 

We use a computer with Intel (R) Xeon (R) 2.40 GHz cpus and 24GB memory, and  MATLAB R2010a. To construct Lasserre's \cite{LAS01}, sparse \cite{WAKI06} and Adaptive SOS problems, we use SparsePOP 2.99 \cite{WAKI08}. To solve the resulting SDP relaxation problems, we use SeDuMi 1.3 \cite{SEDUMI} and SDPT3 4.0 \cite{SDPT3} with the default parameters. The default tolerances for stopping criterion of SeDuMi and SDPT3 are 1.0e-9 and 1.0e-8, respectively.  

To determine whether the optimal value of an  SDP relaxation problem is the exact optimal value of a given POP or not, we use the following two criteria $\epsilon_{\mbox{obj}}$ and $\epsilon_{\mbox{feas}}$: Let $\hat{x}$ be a candidate of an optimal solution of the POP obtained from the  SDP relaxations. We apply a projection of the dual solution of the SDP relaxation problem onto $\mathbb{R}^n$ for obtaining $\hat{x}$ in this section. See \cite{WAKI06} for the details. We define:
\def\eObj{\epsilon_{\mbox{obj}}}
\def\eFeas{\epsilon_{\mbox{feas}}}
\begin{eqnarray*}
\eObj   &:=&\frac{|\mbox{the optimal value of the SDP relaxation} - f(\hat{x})|}{\max\{1, |f(\hat{x})|\}}, \\
\eFeas &:=&\min_{k=1, \ldots, m}\{f_k(\hat{x})\}.
\end{eqnarray*}
If $\eFeas\ge 0$, then $\hat{x}$ is feasible for the POP.
In addition, if $\eObj=0$, then $\hat{x}$ is an optimal solution of the POP and $f(\hat{x})$ is the optimal value of the POP. 

We introduce the following value to indicate the closeness between the obtained values of Lasserre's, sparse and Adaptive SOS relaxations. 
\begin{equation}\label{ratio}
\mbox{Ratio} := \frac{(\mbox{obj. val. of Lasserre's or sparse SDP relax. })}{(\mbox{obj. val. of Adaptive SOS relax.})}=\frac{\rho^*_r}{\rho^*(r)}. 
\end{equation}
If the signs of both optimal values are the same and $\mbox{Ratio}$ is sufficiently close to 1, then the optimal value of  Adaptive SOS relaxation is close to the optimal value of  Lasserre's \redch{}{and sparse SDP relaxations.}  In general, this value is meaningless for measuring the closeness if those signs are different or either of values is zero. Fortunately, those values are not zero and those signs are the same in all numerical experiments in this section.

To reduce the size of the resulting SDP relaxation problems, SparsePOP has functions based on the methods proposed in \cite{KOJ03a, WAKI11b}. 
These methods are closely related to a facial reduction algorithm proposed by Borwein and Wolkowicz \cite{borweinJAMS1,borweinJMAA}, and thus we can expect the numerical stability of the primal-dual interior-point methods for the SDP relaxations may be improved. 
In this section, except for Subsection \ref{subsec:nonarchi}, we apply the method proposed in \cite{WAKI11b}.  

For POPs which have lower and upper bounds on variables, we can strengthen the SDP relaxations by adding valid inequalities based on these bound constraints. In this section, we add them as in \cite{WAKI06}. See Subsection 5.5 in \cite{WAKI06} for the details. 

Table \ref{notation} shows the notation used in the description of numerical experiments in the following subsections.

\begin{table}[htdp]
\caption{Notation}

\begin{center}
{\begin{tabular}{cp{5in}}
\hline
iter. & the number of iterations in SeDuMi and SDPT3\\
rowA, colA & the size of coefficient matrix $A$ in the SeDuMi input format\\
nnzA & the number of nonzero elements in coefficient matrix $A$ in the SeDuMi input format\\
SDPobj & the objective value obtained by SeDuMi for the resulting SDP relaxation problem\\
POPobj & the value of $f$ at a solution $\hat{x}$ retrieved by SparsePOP\\
\#solved &the number of the POPs which are solved by SDP relaxation in 30 problems. If both $\eObj$ and $\eFeas$ are smaller than 1.0e-7, we regard that the SDP relaxation attains the optimal value of the POP. \\
minRatio & minimum value of Ratio defined in (\ref{ratio}) in 30 problems \\
aveRatio & average of Ratio defined in (\ref{ratio}) in 30 problems \\
maxRatio & maximum value of Ratio defined in (\ref{ratio}) in 30 problems \\
$\sec$ & cpu time consumed by SeDuMi or SDPT3 in seconds\\
min.t & minimum cpu time consumed by SeDuMi or SDPT3 in seconds among 30 resulting SDP relaxations\\
ave.t & average cpu time consumed by SeDuMi or SDPT3   in seconds among 30 resulting SDP relaxations\\
max.t & maximum cpu time consumed by SeDuMi or SDPT3   in seconds among 30 resulting SDP relaxations\\
\hline
\end{tabular}}
\end{center}
\label{notation}
\end{table}%

\subsection{Numerical results for POP whose quadratic module is non-archimedean}\label{subsec:nonarchi}

In this subsection, we give the following POP and apply Adaptive SOS relaxation:
\begin{equation}
\label{prestel-delzell}
\inf_{x, y\in\mathbb{R}}\left\{
-x-y: 
\begin{array}{l}
f_1(x, y) := x-0.5\ge 0,\\
f_2(x, y) := y-0.5\ge 0,\\
f_3(x, y) := 0.5-xy\ge 0
\end{array}
\right\}. 
\end{equation}
The optimal value is $-1.5$ and the solutions are $(0.5, 1)$ and $(1, 0.5)$. It was proved in \cite{prestel, WAKI11} that the quadratic module associated with POP (\ref{prestel-delzell}) is non-\redch{A}{a}rchimedean and that all the resulting SDP relaxation problems are weakly infeasible. However, the convergence of computed values of Lasserre's SDP relaxation for POP (\ref{prestel-delzell}) was observed in \cite{WAKI11}. 

In \cite{WAKI11}, it was shown that Lasserre's SDP relaxation \eqref{eq:SOS1}
for \eqref{prestel-delzell} is weakly infeasible.
Since Adaptive SOS relaxation for \eqref{prestel-delzell} has less 
monomials for representing $\sigma_j$'s than that of Lasserre's, the 
resulting SDP relaxation problems are necessarily infeasible.

However, we expect from Thorem \ref{th:conv} that 
Adaptive SOS relaxation attains the optimal value $-1.5$.
Table \ref{subsec5.1_table} provides numerical results 
for Adaptive SOS relaxation based on \eqref{sosProb}. 
In fact, we observe from Table \ref{subsec5.1_table} that 
$\rho^\ast(r)$ obtained by SeDuMi is equal to $-1.5$ at 
$r=7,8,9,10$.
By SDPT3, we observe similar results.

\begin{table}
\caption{The approximate optimal value, cpu time, the number of iterations by SeDuMi and SDPT3}
{\footnotesize
\begin{center}
\begin{tabular}{|cc|c|cc|}
\hline
$r$ & Software & iter.& SDPobj & [$\sec$]\\
\hline
1&SeDuMi&46&-5.9100801e+07&0.31\\ 
&SDPT3&37&-1.8924840e+06&0.57\\\hline
2&SeDuMi&38&-6.8951407e+02&0.29\\ 
&SDPT3&72&-1.1676106e+04&1.28\\\hline
3&SeDuMi&32&-4.2408507e+01&0.22\\ 
&SDPT3&77&-2.0928888e+00&1.43\\\hline
4&SeDuMi&35&-1.2522887e+01&0.30\\ 
&SDPT3&76&-1.8195861e+00&1.74\\\hline
5&SeDuMi&32&-3.5032311e+00&0.39\\ 
&SDPT3&86&-1.6015287e+00&2.65\\\hline
6&SeDuMi&33&-1.8717460e+00&0.48\\ 
&SDPT3&86&-1.5025613e+00&3.43\\\hline
7&SeDuMi&17&-1.5000064e+00&0.47\\ 
&SDPT3&21&-1.5000022e+00&1.18\\\hline
8&SeDuMi&16&-1.5000030e+00&0.58\\ 
&SDPT3&25&-1.5000001e+00&2.03\\\hline
9&SeDuMi&15&-1.5000023e+00&0.75\\ 
&SDPT3&21&-1.4999912e+00&1.95\\\hline
10&SeDuMi&15&-1.5000015e+00&0.99\\ 
&SDPT3&17&-1.5003641e+00&1.89\\\hline
\end{tabular}
\end{center}
}
\label{subsec5.1_table}
\end{table}

\subsection{The difference between Lasserre's and Adaptive SOS  relaxations}\label{subsec:diff}

In this subsection, we show a POP where Adaptive SOS relaxation converges to the optimal value strictly slower than Lasserre's, practically. This POP is available at \cite{GLOBALLIB}, whose name is ``st\_e08.gms". 
\begin{equation}\label{st_e08}
\inf_{x, y\in\mathbb{R}}\left\{
2x + y:
\begin{array}{ll}
f_1(x, y) := xy -1/16\ge 0, & f_2(x, y):= x^2 + y^2-1/4 \ge 0,\\
f_3(x, y) := x \ge 0, & f_4(x, y) :=1-x\ge 0, \\
f_5(x, y) := y \ge 0, & f_6(x, y) :=1-y\ge 0.
\end{array}
\right\}. 
\end{equation}
The optimal value is $(3\sqrt{6} -\sqrt{2})/8 \approx 0.741781958247055$ and solution is $(x^*, y^*) = ((\sqrt{6} -\sqrt{2})/8, (\sqrt{6} +\sqrt{2})/8)$. 
\begin{table}[htdp]
\caption{Numerical results on SDP relaxation problems in Subsection 3.2 by SeDuMi and SDPT3}
{\footnotesize
\begin{center}
\begin{tabular}{|c|c|l|l|}
\hline
\multicolumn{2}{|c|}{}&Lasserre&Adaptive SOS\\
\hline
$r$ & Software & (SDPobj, POPobj$|$ $\epsilon_{\mbox{obj}}, \epsilon_{\mbox{feas}}$ $|$ [$\sec$])& (SDPobj, POPobj$|$ $\epsilon_{\mbox{obj}}, \epsilon_{\mbox{feas}}$ $|$ [$\sec$])  \\
\hline
1 & SeDuMi &(0.00000e+00, 0.00000e+00$|$ 0.0e+00, -1.0e+00$|$ 0.02) &(0.00000e+00, 0.00000e+00$|$ 0.0e+00, -1.0e+00$|$ 0.02 )  \\
   & SDPT3  &(-1.16657e-09, 5.89142e-10$|$ 1.8e-09, -1.0e+00$|$ 0.14)&(-1.16657e-09,  5.89142e-10$|$  1.8e-09, -1.0e+00$|$  0.06) \\
\hline
2 & SeDuMi &(3.12500e-01, 3.12500e-01$|$ -9.5e-10, -8.4e-01$|$  0.09) &(2.69356e-01, 2.69356e-01$|$ -1.7e-10, -9.3e-01$|$  0.09)  \\
   &SDPT3 &(3.12500e-01, 3.12500e-01$|$  2.0e-09 ,  -8.4e-01$|$  0.22)&(2.69356e-01,  2.69356e-01$|$  1.1e-09, -9.3e-01$|$  0.21) \\
\hline
3 & SeDuMi &(7.41782e-01, 7.41782e-01$|$ -2.0e-11, -1.1e-09$|$   0.15) &(3.06312e-01, 3.06312e-01$|$ -1.1e-09, -8.3e-01$|$   0.13)  \\
   & SDPT3 &(7.41782e-01,  7.41782e-01$|$  2.0e-08,  0.0e+00$|$  0.26)&(3.06312e-01,  3.06312e-01$|$  4.6e-09,  -8.3e-01$|$  0.25) \\
\hline
4 & SeDuMi &(7.41782e-01, 7.41782e-01$|$ 1.1e-10, -1.5e-09$|$   0.15) &(7.29855e-01, 7.29855e-01$|$ -1.2e-07, -4.9e-02$|$   0.24)  \\
   &SDPT3 &(7.41782e-01,  7.41782e-01$|$  2.8e-09,  0.0e+00$|$  0.34)&(7.29855e-01, 7.29855e-01$|$  2.5e-08, -4.9e-02$|$  0.36) \\
\hline
5 & SeDuMi &(7.41782e-01, 7.41782e-01$|$ 8.3e-11, -4.5e-10$|$   0.19) &(7.36195e-01, 7.36194e-01$|$ -9.5e-07, -4.2e-02$|$   0.33) \\
   &SDPT3 &(7.41782e-01,  7.41782e-01$|$  -6.3e-10,  0.0e+00$|$  0.72) &(7.36195e-01,  7.36195e-01$|$  5.3e-08,  -4.2e-02$|$  0.50) \\
\hline
6 & SeDuMi &(7.41782e-01, 7.41782e-01$|$ 2.3e-11, -6.1e-11$|$   0.27) &(7.41782e-01, 7.41782e-01$|$ -1.0e-09, -6.6e-09$|$   0.20) \\
   &SDPT3 &(7.41782e-01,  7.41782e-01$|$  3.4e-10,  0.0e+00$|$  1.02) &(7.41782e-01,  7.41782e-01$|$  -4.7e-11,  0.0e+00$|$  0.98) \\
\hline
\end{tabular}
\end{center}
}
\label{subsec5.2_info_table}
\end{table}

Table \ref{subsec5.2_info_table} 
show the numerical results of SDP relaxations for POP (\ref{st_e08}) by SeDuMi and SDPT3. We observe that Lasserre's SDP relaxation attains the optimal value of (\ref{st_e08}) by relaxation order $r=3$, while Adaptive SOS relaxation attains it only at the relaxation order by $r=6$. 

\subsection{Numerical results for detecting the copositivity}
\label{subsec:cop}

The symmetric matrix $A$ is said to be {\it copositive} if $x^TAx \ge 0$ for all $x\in\mathbb{R}^n_+$. We can formulate the problem for detecting whether a given matrix is copositive, as follows:
\begin{equation}\label{cop}
\inf_{x\in\mathbb{R}^n}\left\{
x^TAx :
f_i(x) := x_i\ge 0 \ (i=1, \ldots, n), f_{n+1}(x) := 1- \sum_{i=1}^n x_i = 0, 
\right\}. 
\end{equation} 
If the optimal value of this problem is nonnegative, then $A$ is copositive. 
\begin{table}
\caption{Information on SDP relaxations problems in Subsection 3.3 by SeDuMi and SDPT3}
{\footnotesize 
\begin{center}
\begin{tabular}{|c|c|l|l||l|}
\hline
\multicolumn{2}{|c|}{}&Lasserre&Adaptive SOS&\\
\hline
$n$ & Software & (\#solved $|$ min.t, ave.t, max.t)&(\#solved $|$ min.t, ave.t, max.t) &(minR, aveR, maxR) \\
\hline
5 & SeDuMi & (30 $|$ 0.14 0.18 0.50) &(30 $|$ 0.12 0.16 0.20)&(1.0,  1.0, 1.0)\\
  & SDPT3 & (30 $|$ 0.40 0.44 0.85) & (30 $|$ 0.34 0.42 0.53) & (1.0, 1.0, 1.0)\\
\hline
10 & SeDuMi & (29 $|$ 0.36 0.42 0.50)&(29 $|$ 0.23 0.31 0.42)& (1.0, 1.0, 1.0)\\
& SDPT3 & (29 $|$ 0.73 1.00 1.48)&(30 $|$ 0.66 0.88 1.23)&  (1.0, 1.0, 1.0)\\
\hline
15 & SeDuMi & (30 $|$ 1.59 1.99 2.52) &(30 $|$ 0.75 0.99 1.31)& (1.0, 1.0, 1.0)\\
& SDPT3 & (29 $|$ 2.91 3.40 4.73)&(23 $|$ 1.58 2.04 2.80)& (1.0,  1.0,  1.0)\\
\hline
20 & SeDuMi & (30 $|$ 10.22 14.06 19.98)&(30 $|$ 4.47 6.02 7.72)& (1.0, 1.0, 1.0)\\
& SDPT3 & (26 $|$ 11.40 16.23 19.73)&(1 $|$ 6.65 8.64 11.32)& (1.0, 1.0, 1.0)\\
\hline
25 & SeDuMi & (29 $|$ 215.94 263.88 336.96)&(29 $|$ 49.69 66.63 84.07)& (1.0, 1.0, 1.0)\\
& SDPT3 & (20 $|$ 51.53 64.31 77.35)&(4 $|$ 26.91 36.06 44.74)& (1.0, 1.0, 1.0)\\
\hline
30 & SeDuMi & (27 $|$ 1970.59 2322.30 2930.30)&(28 $|$ 1031.91 1198.05 1527.01)& (1.0, 1.0, 1.0)\\
& SDPT3 & (0 $|$ 136.59 401.23 1184.76)&(0 $|$ 92.96 165.22 295.23)& (0.4, 1.0, 1.6)\\
\hline
\end{tabular}
\end{center}
}
\label{subsec5.3_table}
\end{table}
In this experiment, we solve 30 problems generated randomly. In particular, the coefficients of all diagonal of $A$ are set to be $\sqrt{n}/2$ and the other coefficients are chosen from [-1, 1] uniformly. In addition, since the positive semidefiniteness implies the copositivity, we chose the matrices $A$ which are not positive semidefinite. 

We apply Lasserre's and Adaptive SOS relaxations with relaxation order $r=2$. Table \ref{subsec5.3_table} shows the numerical results by SeDuMi and SDPT3 for (\ref{cop}), respectively.  We observe the following. 
\begin{itemize}
\item SDPT3 fails to solve almost all problems (\ref{cop}), while SeDuMi solves them for $n=20, 25, 30$. 
In particular,  Adaptive SOS relaxations return the optimal values of the original problems although it is no stronger than Lasserre's theoretically. 
\item SeDuMi solves Adaptive SOS relaxation problems faster than Lasserre's because the sizes of Adaptive SOS relaxation problems are smaller than those of Lasserre's. 
\item SDPT3 cannot solve any problems with $n=30$ by Lasserre's and Adaptive SOS relaxation although it terminates faster than SeDuMi. In particular, for almost all SDP relaxation problems, SDPT3 returns the message ``stop: progress is bad" or ``stop: progress is slow" and terminates. This means that it is difficult for SDPT3 to solve those SDP relaxation problems numerically. 
\end{itemize}

\subsection{Numerical results for BoxQP}
\label{subsec:boxqp}
In this subsection, we solve BoxQP:
\begin{equation}\label{boxqp}
\inf_{x\in\mathbb{R}^n}\left\{
x^TQx + c^Tx :
0\le x_i\le 1 \ (i=1, \ldots, n)
\right\}, 
\end{equation}
where each element in $Q\in\mathbb{S}^n$ and $c\in\mathbb{R}^n$ is chosen from [-50, 50] uniformly.  In particular, we vary the number $n$ of the variables in (\ref{boxqp}) and the density of $Q, c$. In this subsection, we compare Adaptive SOS relaxation based on Theorem \ref{th:sparsemain} with sparse SDP relaxation \cite{WAKI06} instead of Lasserre's. Indeed, when the density of $Q$  is small, the BoxQP has sparse structure, and thus sparse SDP relaxation is more effective than Lasserre's. 

\begin{table}
\caption{Information on SDP relaxation problems in Subsection 3.4 with density 0.2 by SeDuMi and SDPT3}
{\footnotesize 
\begin{center}
\begin{tabular}{|c|c|l|l||l|}
\hline
\multicolumn{2}{|c|}{}&Sparse&Adaptive SOS&\\
\hline
$n$ & Software & (\#solved $|$ min.t, ave.t, max.t)&(\#solved $|$ min.t, ave.t, max.t) &(minR, aveR, maxR) \\
\hline
 5 & SeDuMi & (23 $|$0.15, 0.24, 0.48)& (23 $|$ 0.14, 0.23, 0.52)&(0.00072, 12.34638, 342.39518)\\
& SDPT3 &(23 $|$0.20, 0.37, 2.48)&(22 $|$ 0.19, 0.26, 0.34)&(0.00072, 0.97463, 1.24265)\\
\hline
10 & SeDuMi & (13 $|$0.33, 0.55, 0.70)& (12 $|$ 0.28, 0.40, 0.53) &(0.97227, 0.99609, 1.00000)\\
& SDPT3 &(12  $|$0.28, 0.51, 0.62)&(12 $|$ 0.22, 0.29, 0.39)&(0.97227, 0.99609, 1.00000)\\
\hline
15 & SeDuMi & (14  $|$0.57, 0.95, 1.68)& ( 3 $|$ 0.42, 0.63, 0.85)&(0.96590, 0.99172, 1.00000)\\
& SDPT3 &(14  $|$0.54, 0.93, 1.22)& ( 3  $|$ 0.43, 0.57, 0.76)&(0.96590, 0.99172, 1.00000)\\
\hline
20 & SeDuMi & (11 $|$1.40, 2.57, 5.32)&( 0  $|$ 0.80, 0.97, 1.27) &(0.94812, 0.98422, 0.99978)\\
& SDPT3 &(10 $|$1.41, 2.31, 3.55) &( 0 $|$ 0.55, 0.69, 1.01)&(0.94812, 0.98422, 0.99978)\\
\hline
25 & SeDuMi & ( 7 $|$2.57, 5.15, 10.03)& ( 0 $|$ 0.95, 1.09, 1.42)&(0.94333, 0.97591, 0.99923)\\
& SDPT3 &( 6 $|$4.60, 7.24, 12.46) &( 0 $|$ 0.59, 0.85, 1.31)&(0.94333, 0.97591, 0.99923)\\
\hline
30 & SeDuMi & (12 $|$3.43, 15.60, 26.86)&( 0  $|$ 1.27, 1.51, 2.02)&(0.93773, 0.97542, 0.99843)\\
& SDPT3 &(10 $|$8.02, 22.87, 38.42) &( 0 $|$ 0.94, 1.33, 1.67)&(0.93773, 0.97542, 0.99843)\\
\hline
35 & SeDuMi & (12 $|$26.57, 67.79, 143.06)&( 0 $|$ 1.77, 2.15, 3.33) &(0.93271, 0.97236, 0.99648)\\
& SDPT3 & ( 9 $|$44.14, 80.48, 135.30)&( 0 $|$ 1.06, 1.83, 2.63)&(0.93271, 0.97236, 0.99648)\\
\hline
40 & SeDuMi & Not solved&( 0  $|$ 2.47, 2.89, 3.57)&(--, --, --)\\
& SDPT3& Not solved&( 0 $|$ 2.13, 3.13, 3.87)  &(--, --, --)\\
\hline
45 & SeDuMi & Not solved& ( 0  $|$ 3.58, 4.17, 5.51) &(--, --, --)\\
& SDPT3 &Not solved&( 0  $|$ 4.12, 5.09, 6.35)&(--, --, --)\\
\hline
50 & SeDuMi &Not solved& ( 0 $|$ 5.30, 7.02, 9.48)&(--, --, --)\\
& SDPT3 &Not solved&( 0  $|$ 5.19, 6.83, 8.34)&(--, --, --)\\
\hline
55 & SeDuMi &Not solved&  ( 0 $|$ 8.75, 10.43, 12.23)&(--, --, --)\\
& SDPT3 &Not solved&( 0  $|$ 8.31, 10.77, 13.60)&(--, --, --)\\
\hline
60 & SeDuMi & Not solved&  ( 0 $|$ 12.21, 15.16, 19.59)&(--, --, --)\\
& SDPT3 &Not solved&( 0 $|$ 12.62, 16.57, 22.44)&(--, --, --)\\
\hline
\end{tabular}
\end{center}
}
\label{subsec5.4_table_density0.2}
\end{table}

\begin{table}
\caption{Information on SDP relaxation problems in Subsection 3.4 with density 0.4 by SeDuMi and SDPT3}
{\footnotesize 
\begin{center}
\begin{tabular}{|c|c|l|l||l|}
\hline
\multicolumn{2}{|c|}{}&Sparse&Adaptive SOS&\\
\hline
$n$ & Software & (\#solved $|$ min.t, ave.t, max.t)&(\#solved $|$ min.t, ave.t, max.t) &(minR, aveR, maxR) \\
\hline
 5 & SeDuMi & (24  $|$0.14, 0.19, 0.28)& (22  $|$ 0.13, 0.17, 0.27)&(0.98678, 0.99849, 1.00000)\\
& SDPT3 &(23  $|$0.23, 0.27, 0.35)&(22 $|$ 0.17, 0.23, 0.34)&(0.98678, 0.99849, 1.00000 )\\
\hline
10 & SeDuMi & (19  $|$0.28, 0.49, 0.77)&( 9  $|$ 0.25, 0.35, 0.46) &(0.95400, 0.98958, 1.00000)\\
& SDPT3 &(18  $|$0.32, 0.53, 0.86) &( 7  $|$ 0.26, 0.33, 0.52)&(0.95400, 0.98958, 1.00000)\\
\hline
15 & SeDuMi & (13  $|$0.76, 1.21, 2.50)&( 3  $|$ 0.46, 0.56, 0.65)&(0.95219, 0.98580, 1.00000)\\
& SDPT3 & (13  $|$0.84, 1.32, 2.26)&( 3  $|$ 0.37, 0.54, 0.81)&(0.95219, 0.98580, 1.00000)\\
\hline
20 & SeDuMi & (11  $|$2.10, 3.51, 5.45)&( 0  $|$ 0.70, 0.79, 0.97) &(0.94457, 0.97953, 0.99933)\\
& SDPT3 & (11  $|$3.22, 5.61, 8.30)&( 0  $|$ 0.50, 0.73, 1.01)&(0.94457, 0.97953, 0.99933)\\
\hline
25 & SeDuMi & (11  $|$6.65, 13.88, 24.32)&( 0  $|$ 1.02, 1.13, 1.28)&(0.92917, 0.96999, 0.99596)\\
& SDPT3 & (10  $|$11.48, 21.00, 30.98)&( 0  $|$ 0.69, 1.03, 1.47)&(0.92917, 0.96999, 0.99596)\\
\hline
30 & SeDuMi & (14  $|$27.25, 60.67, 108.22)&( 0  $|$ 1.31, 1.62, 2.26)&(0.92761, 0.97283, 0.99608)\\
& SDPT3 & (12  $|$43.33, 66.25, 95.80)&( 0  $|$ 1.29, 1.71, 2.22)&(0.92761, 0.97283, 0.99608)\\
\hline
35 & SeDuMi & ( 8  $|$76.07, 328.08, 589.43)&( 0  $|$ 2.11, 2.42, 2.95)&(0.93669, 0.96707, 0.99717)\\
& SDPT3 &( 6  $|$116.23, 218.61, 322.82) &( 0  $|$ 2.21, 2.87, 5.03)&(0.93669, 0.96707, 0.99717)\\
\hline
40 & SeDuMi & Not solved&( 0  $|$3.11, 3.54, 4.69) &(--, --, --)\\
& SDPT3 &Not solved &( 0  $|$3.29, 4.50, 5.39)&(--, --, --)\\
\hline
45 & SeDuMi & Not solved&( 0  $|$4.99, 5.79, 7.10) &(--, --, --)\\
& SDPT3 &Not solved & ( 0  $|$5.43, 6.89, 8.85)&(--, --, --)\\
\hline
50 & SeDuMi & Not solved&( 0  $|$7.09, 8.47, 11.58)&(--, --, --)\\
& SDPT3 & Not solved &( 0  $|$9.09, 11.30, 15.02)&(--, --, --)\\
\hline
55 & SeDuMi & Not solved&( 0  $|$11.84, 14.34, 17.72) &(--, --, --)\\
& SDPT3 &Not solved &( 0  $|$14.09, 18.30, 22.13)&(--, --, --)\\
\hline
60 & SeDuMi & Not solved&( 0  $|$19.33, 24.23, 29.13) &(--, --, --)\\
& SDPT3 &Not solved &( 0  $|$19.45, 22.96, 26.65)&(--, --, --)\\
\hline
\end{tabular}
\end{center}
}
\label{subsec5.4_table_density0.4}
\end{table}

\begin{table}
\caption{Information on SDP relaxation problems in Subsection 3.4 with density 0.6 by SeDuMi and SDPT3}
{\footnotesize 
\begin{center}
\begin{tabular}{|c|c|l|l||l|}
\hline
\multicolumn{2}{|c|}{}&Sparse&Adaptive SOS&\\
\hline
$n$ & Software & (\#solved $|$ min.t, ave.t, max.t)&(\#solved $|$ min.t, ave.t, max.t) &(minR, aveR, maxR) \\
\hline
 5 & SeDuMi & (27 $|$0.13, 0.22, 0.54)&(25 $|$ 0.12, 0.17, 0.38)&(0.93673, 0.99543, 1.00000)\\
& SDPT3 &(26 $|$0.21, 0.26, 0.33) &(25 $|$ 0.18, 0.21, 0.29)&(0.93673, 0.99543, 1.00000)\\
\hline
10 & SeDuMi & (19 $|$0.36, 0.68, 1.26)&( 6 $|$ 0.33, 0.48, 0.79) &(0.94709, 0.98678, 1.00000)\\
& SDPT3 & (18 $|$0.37, 0.48, 0.72)&( 6 $|$  0.25, 0.31, 0.40)&(0.94709, 0.98678, 1.00000)\\
\hline
15 & SeDuMi & (14 $|$0.71, 1.52, 3.70)&( 6 $|$ 0.42, 0.61, 1.01) &(0.95463, 0.98581, 1.00000)\\
& SDPT3 &(14 $|$0.77, 1.33, 2.04) &( 6 $|$ 0.34, 0.41, 0.51)&(0.95463, 0.98581, 1.00000)\\
\hline
20 & SeDuMi & (13 $|$1.92, 5.18, 7.99)&( 2 $|$ 0.72, 0.91, 1.56)&(0.92378, 0.97521, 1.00000)\\
& SDPT3 &(11 $|$2.25, 5.54, 8.21) &( 2 $|$  0.52, 0.61, 0.75)&(0.92378, 0.97521, 1.00000)\\
\hline
25 & SeDuMi & (15 $|$9.56, 29.31, 57.08)&( 0 $|$ 1.03, 1.24, 1.94) &(0.92768, 0.96827, 0.99715)\\
& SDPT3 &(12 $|$15.55, 26.06, 40.61) &( 0 $|$ 0.75, 0.93, 1.19)&(0.92768, 0.96827, 0.99715)\\
\hline
30 & SeDuMi & (11 $|$50.72, 168.53, 368.04)&( 0 $|$ 1.56, 1.97, 2.99) &(0.93048, 0.96888, 0.99470)\\
& SDPT3 &( 9 $|$42.25, 90.31, 140.94) &( 0 $|$ 1.27, 1.50, 2.10)&(0.93048, 0.96888, 0.99470)\\
\hline
35 & SeDuMi & (12 $|$510.67, 964.20, 1489.56)&( 0 $|$ 2.52, 3.11, 4.27)&(0.90892, 0.95875, 0.99301)\\
& SDPT3 &(11 $|$217.87, 303.90, 366.57) &( 0 $|$ 2.16, 2.55, 3.09)&(0.90892, 0.95875, 0.99301)\\
\hline
40 & SeDuMi & Not solved&( 0 $|$3.77, 4.34, 5.77)&(--, --, --)\\
& SDPT3 &Not solved &( 0 $|$3.37, 4.24, 5.12)&(--, --, --)\\
\hline
45 & SeDuMi & Not solved&( 0 $|$6.08, 6.91, 8.33)&(--, --, --)\\
& SDPT3 &Not solved &( 0 $|$5.63, 7.07, 9.33)&(--, --, --)\\
\hline
50 & SeDuMi & Not solved&( 0 $|$8.97, 10.66, 12.82)&(--, --, --)\\
& SDPT3 &Not solved &( 0 $|$8.87, 10.59, 11.84)&(--, --, --)\\
\hline
55 & SeDuMi & Not solved&( 0 $|$13.95, 17.13, 20.71)&(--, --, --)\\
& SDPT3 &Not solved &( 0 $|$10.26, 13.64, 20.92)&(--, --, --)\\
\hline
60 & SeDuMi & Not solved&( 0 $|$21.94, 25.42, 30.36)&(--, --, --)\\
& SDPT3 & Not solved&( 0 $|$15.48, 19.66, 27.26)&(--, --, --)\\
\hline
\end{tabular}
\end{center}
}
\label{subsec5.4_table_density0.6}
\end{table}

\begin{table}
\caption{Information on SDP relaxation problems in Subsection 3.4 with density 0.8 by SeDuMi and SDPT3}
{\footnotesize 
\begin{center}
\begin{tabular}{|c|c|l|l||l|}
\hline
\multicolumn{2}{|c|}{}&Sparse&Adaptive SOS&\\
\hline
$n$ & Software & (\#solved $|$ min.t, ave.t, max.t)&(\#solved $|$ min.t, ave.t, max.t) &(minR, aveR, maxR) \\
\hline
 5 & SeDuMi & (25 $|$0.15, 0.19, 0.34)&(22 $|$  0.13, 0.17, 0.24)&(0.94896, 0.99548, 1.00000)\\
& SDPT3 &(25 $|$0.22, 0.27, 0.37) &(22 $|$ 0.18, 0.22, 0.29)&(0.94896, 0.99548, 1.00000)\\
\hline
10 & SeDuMi & (20 $|$0.36, 0.54, 0.80)&(11 $|$  0.26, 0.37, 0.52)&(0.96388, 0.99365, 1.00000)\\
& SDPT3 & (20 $|$0.40, 0.62, 0.99)&(10 $|$ 0.29, 0.39, 0.59)&(0.96388, 0.99365, 1.00000)\\
\hline
15 & SeDuMi & (14 $|$0.93, 1.67, 2.93)&( 1 $|$  0.50, 0.59, 0.71)&(0.94514, 0.98537, 1.00000)\\
& SDPT3 &(12 $|$1.29, 1.85, 2.63) &( 1 $|$  0.42, 0.51, 0.71)&(0.94514, 0.98537, 1.00000)\\
\hline
20 & SeDuMi & (14 $|$2.51, 5.22, 8.98)&( 2 $|$ 0.66, 0.85, 1.15)&(0.95261, 0.98061, 1.00000)\\
& SDPT3 &(12 $|$4.50, 6.70, 9.35) &( 2  $|$ 0.56, 0.76, 1.13)&(0.95261, 0.98061, 1.00000)\\
\hline
25 & SeDuMi & (10 $|$10.64, 23.57, 56.02)&( 0  $|$ 1.13, 1.25, 1.52)&(0.95060, 0.97500, 0.99997)\\
& SDPT3 & (10 $|$14.13, 26.81, 44.75)&( 0 $|$ 0.87, 1.11, 1.66)&(0.95060, 0.97500, 0.99997)\\
\hline
30 & SeDuMi & (11 $|$42.70, 156.60, 507.20)&( 0 $|$ 1.68, 1.89, 2.18)&(0.94199, 0.96738, 0.99484)\\
& SDPT3 & ( 9 $|$53.52, 104.12, 173.49)&( 0 $|$ 1.43, 1.88, 2.49)&(0.94199, 0.96738, 0.99484 )\\
\hline
35 & SeDuMi & (15 $|$185.51, 1000.24, 2158.08)&( 0 $|$  2.66, 2.89, 3.15)&(0.92313, 0.96254, 0.99485)\\
& SDPT3 &(12 $|$157.31, 337.69, 508.43) &( 0  $|$ 2.52, 2.99, 3.60)&(0.92313, 0.96258, 0.99485)\\
\hline
40 & SeDuMi & Not solved& ( 0 $|$4.45, 4.89, 6.34)&(--, --, --)\\
& SDPT3 &Not solved & ( 0 $|$4.11, 5.22, 6.66)&(--, --, --)\\
\hline
45 & SeDuMi & Not solved&( 0 $|$6.52, 7.63, 8.86)&(--, --, --)\\
& SDPT3 & Not solved &( 0 $|$7.00, 8.05, 9.51)&(--, --, --)\\
\hline
50 & SeDuMi & Not solved&( 0 $|$10.45, 11.70, 13.89) &(--, --, --)\\
& SDPT3 &Not solved &( 0 $|$10.57, 12.65, 15.41)&(--, --, --)\\
\hline
55 & SeDuMi & Not solved&( 0 $|$15.96, 19.55, 24.40) &(--, --, --)\\
& SDPT3 &Not solved &( 0 $|$11.84, 16.07, 21.26)&(--, --, --)\\
\hline
60 & SeDuMi & Not solved&( 0 $|$26.31, 32.04, 36.89) &(--, --, --)\\
& SDPT3 &Not solved &( 0 $|$17.69, 22.33, 27.93)&(--, --, --)\\
\hline
70 & SeDuMi & Not solved&( 0 $|$69.62, 91.01, 123.14) &(--, --, --)\\
& SDPT3 & Not solved&( 0 $|$26.30, 34.00, 45.75)&(--, --, --)\\
\hline
80 & SeDuMi & Not solved&( 0 $|$182.40, 218.82, 268.42) &(--, --, --)\\
& SDPT3 &Not solved &( 0 $|$46.87, 52.48, 59.51)&(--, --, --)\\
\hline
90 & SeDuMi & Not solved&( 0 $|$406.85, 478.44, 619.49) &(--, --, --)\\
& SDPT3 & Not solved&( 0 $|$77.36, 91.34, 107.29)&(--, --, --)\\
\hline
100 & SeDuMi & Not solved&( 0 $|$844.15, 943.74, 1138.27) &(--, --, --)\\
& SDPT3 &Not solved &( 0 $|$130.50, 148.36, 172.25)&(--, --, --)\\
\hline
\end{tabular}
\end{center}
}
\label{subsec5.4_table_density0.8}
\end{table}

We observe the following from Table  \ref{subsec5.4_table_density0.2}. 
\begin{itemize}
\item Sparse SDP relaxation obtains the optimal solution for some BoxQPs, while Adaptive SOS relaxation cannot. 
\item Adaptive SOS relaxation solves the resulting SDP problems approximately 10 $\sim$ 30 times faster than Lasserre's. 
\item The values obtained by Adaptive SOS relaxation  are within 10\% of Sparse SDP relaxation, except for $n=5$. 
\end{itemize}

\subsection{Numerical results for Bilinear matrix inequality eigenvalue problems}\label{subsec:bmi}
In this subsection, we solve the binary matrix inequality eigenvalue problems. 
\begin{equation}\label{BMIEP}
 \inf_{s\in\mathbb{R}, x\in\mathbb{R}^n, y\in\mathbb{R}^m} \SET{s}{sI_k-B_k(x,y) \in \S_+^k, x\in [0,1]^n, y\in[0,1]^m}, 
\end{equation}
where we define for $k\in\N$, $x\in\R^n$ and $y\in\R^m$: 
\[
 B_k(x, y) = \sum_{i=1}^n\sum_{j=1}^m B_{ij} x_i y_j
   + \sum_{i=1}^n B_{i0} x_i 
   + \sum_{j=1}^m B_{0j} y_j 
   +B_{00},
\]
where $B_{ij} (i=0,\ldots,n, j=0,\ldots,m)$ are
$k\times k$ symmetric matrices. In this numerical experiment, each element of  $B_{ij}$ is chosen from $[-1, 1]$ uniformly. (\ref{BMIEP}) is the problem of minimizing the maximum eigenvalue of $B_k(x,y)$
keeping $B_k(x,y)$ positive semidefinite. 

We apply Lasserre and Adaptive SOS relaxations with relaxation order $r=3$. Tables \ref{subsec5.5_table} shows the numerical results for BMIEP (\ref{BMIEP}) with $k=5, 10$ by SeDuMi and SDPT3, respectively. 

\begin{table}
\caption{Information on SDP relaxation problems in Subsection \ref{subsec:bmi} by SeDuMi and SDPT3}
{\footnotesize 
\begin{center}
\begin{tabular}{|c|c|l|l||l|}
\hline
\multicolumn{2}{|c|}{}&Lasserre&Adaptive SOS&\\
\hline
$(n ,m, k)$  & Software & (\#solved $|$ min.t, ave.t, max.t)&(\#solved $|$ min.t, ave.t, max.t) &(minR, aveR, maxR) \\
\hline
(1, 1,  5) & SeDuMi & (21 $|$ 0.11, 0.16, 0.25)&(16  $|$ 0.10, 0.17, 0.29)&(1.00000, 1.00103, 1.02352)\\
& SDPT3 &(21 $|$ 0.28, 0.36, 0.42)&(16 $|$ 0.26, 0.32, 0.41)&(1.00000, 1.00103,  1.02352)\\
\hline
(1, 1, 10) & SeDuMi &(20 $|$ 0.12, 0.16, 0.21)&(18  $|$ 0.11, 0.16,  0.30)&(1.00000, 1.00018, 1.00450)\\
& SDPT3 &(20 $|$ 0.32, 0.37, 0.47) &(18 $|$ 0.26, 0.33, 0.44)&(1.00000,  1.00018,  1.00450)\\
\hline
(3, 3,  5) & SeDuMi & (3 $|$ 1.95,  3.49,  4.74) &(1  $|$ 0.43, 0.63,  0.81)&(0.878394,  1.01520,  1.20254) \\
& SDPT3 & (3 $|$ 4.81, 7.12, 8.77)&(1 $|$ 0.67, 0.97, 1.14)&(0.878394,  1.01520,  1.20254 )\\
\hline
(3, 3, 10) & SeDuMi & (0 $|$ 2.46, 3.89, 4.77)&(0  $|$ 0.54, 0.69, 0.93)&(1.00000,  1.00407,  1.01243)\\
& SDPT3 &(0  $|$ 5.51, 7.63, 8.95) &(0 $|$ 0.88, 1.04, 1.16)&(1.00000,  1.00407,  1.01243 )\\
\hline
(5, 5,  5) &SeDuMi  & (0 $|$ 219.93, 350.02, 545.81)&(0  $|$ 8.25, 10.99, 14.08)& (0.649823, 1.04081, 1.26310)\\
& SDPT3 & (0 $|$ 160.89,  247.24, 298.97)&(0 $|$ 4.45, 5.50, 6.97)&(0.649823,  1.04081, 1.26310)\\
\hline
(5, 5, 10) & SeDuMi & (0 $|$ 285.21, 420.27, 509.31)&(0 $|$  7.96,  10.53, 15.04)&(1.00000, 1.01445, 1.02818)\\
& SDPT3 & (0 $|$ 217.48, 276.67, 309.27)&(0 $|$ 4.34, 5.37, 6.66)&(1.00000,  1.01445,   1.02818)\\
\hline
\end{tabular}
\end{center}
}
\label{subsec5.5_table}
\end{table}



%

We observe the following:
\begin{itemize}
\item SDPT3 solves SDP relaxation problems faster than SeDuMi for $(n, m) = (5, 5)$. 
\item  Adaptive SOS relaxation can solve the resulting SDP problems faster than Lasserre's. In particular, SDPT3 works efficiently for Adaptive SOS relaxation for BMIEP (\ref{BMIEP}). 

\end{itemize}

\section{Extensions}\label{sec:extensions}
In this section, we give three extensions of Theorem \ref{th:main} and present some related work to Theorem \ref{th:main}.  

\subsection{Sums of squares of rational polynomials}\label{subsec:rational}
We can extend part i. of Theorem \ref{th:main} with sums of squares of rational polynomials. We assume that for all $j=1, \ldots, m$, there exists  $g_j\in\mathbb{R}[x]$ such that $|f_j(x)|\le g_j(x)$ and $g_j(x)\neq 0$ for all $x\in B$. We define
\[
\tilde{\psi}_r(x) = - \sum_{j=1}^m f_j(x) \left(1-\frac{f_j(x)}{g_j(x)}\right)^{2r}
\]
for all $r\in\mathbb{N}$. Then, we can prove the following corollary by using almost the same arguments as Theorem \ref{th:main}. 
\begin{corollary} \label{coro}
Suppose that for $\rho\in \R$, $f(x) -\rho > 0$ for every $x\in \bar{K}$,
i.e., $\rho$ is an lower bound of $f^\ast$. Then there exists $\tilde{r}\in\mathbb{N}$ such that for all $r\ge\tilde{r}$, $f-\rho +  \tilde{\psi}_r$ is positive over $B$. 
\end{corollary}

It is difficult to apply Corollary \ref{coro} to the framework of SDP relaxations, because we deal with rational polynomials in $\tilde{\psi}_r$. However, we may be able to reduce the degrees of sums of squares in $\tilde{\psi}_r$ by using Corollary \ref{coro}. For instance, we consider $f_1(x) = 1-x^4$ and $B=[-1,1]$. Choose $g_1(x) = 2(1+x^2)$.  Then $g_1$ dominates $|f_1|$ over $B$, {\it i.e.}, $|f_1(x)|\le g_1(x)$ for all $x\in B$. We have
\[
\tilde{\psi}_r(x) = -(1-x^4)\left(
1 - \frac{1-x^4}{2(1+x^2)}
\right)^{2r} = -(1-x^4)\left(
1 - \frac{1-x^2}{2}
\right)^{2r}, 
\]
and the degree of $\tilde{\psi}$ in Corollary \ref{coro} is $4r$, while the degree of $\psi$  in Theorem \ref{th:main} is $8r$. 

\subsection{Extension to POP with correlative sparsity}\label{subsec:csp}
In \cite{WAKI06}, the authors introduced the notion of correlative sparsity for POP (\ref{eq:POP0}), and  proposed a  sparse SDP relaxation that exploits the correlative sparsity. They then demonstrated that the sparse SDP relaxation outperforms Lasserre's SDP relaxation. The sparse SDP relaxation is implemented in \cite{WAKI08} and its source code is freely available. 

We give {some of} the definition of the correlative sparsity for POP (\ref{eq:POP0}). For this, we use {an} $n\times n$ symbolic symmetric matrix $R$, whose elements are either $0$ or $\star$ representing a nonzero value. We assign either $0$ or $\star$ as follows:
\[
R_{k, \ell} =\left\{
\begin{array}{ll}
\star &\mbox{if } k=\ell, \\
 \star &\mbox{if } \alpha_k\ge 1 \mbox{ and } \alpha_{\ell}\ge 1 \mbox{ for some }\alpha\in \mathcal{F}, \\
 \star &\mbox{if } x_k \mbox{ and } x_{\ell} \mbox{ are involved in the polynomial } f_j  \mbox{ for some } j=1, \ldots, m, \\
0 &\mbox{o.w.}
\end{array}
\right.
\]
POP (\ref{eq:POP0}) is said to be \textit{correlatively sparse} if the matrix $R$ is sparse. 

We give some of the details of the sparse SDP relaxation proposed in \cite{WAKI06} for the sake of completeness. We construct an undirected graph $G=(V, E)$ from $R$. Here $V:=\{1, \ldots, n\}$ and $E:=\{(k, \ell) :  R_{k, \ell}=\star\}$.  After applying the chordal extension to $G=(V, E)$, 
 we generate all maximal cliques $C_1, \ldots, C_p$ of the extension $G=(V, \tilde{E})$ with $E\subseteq\tilde{E}$.  See \cite{FUKUDA00, WAKI06} and references therein for the details of the construction of the chordal extension.   For a finite set $C\subseteq\mathbb{N}$, $x_C$ denotes the subvector which consists of $x_i \ (i\in C)$.  For all $f_1, \ldots, f_m$ in POP (\ref{eq:POP0}), $F_j$ denotes the set of indices whose variables are involved in $f_j$, {\it i.e.}, $F_j := \{
  i\in\{1, \ldots, n\}: \alpha_i\ge 1 \ \mbox{for some } \alpha\in\mathcal{F}_j\}$. For a finite set $C\subseteq\mathbb{N}$,  the sets  $\Sigma_{r, C}$ and $\Sigma_{\infty, C}$ denote the subsets of $\Sigma_r$ as follows:
  \begin{eqnarray*}
  \Sigma_{r, C} &:=&\left\{\sum_{k=1}^{q}g_k(x)^2: \forall k=1, \ldots, q, g_k\in\mathbb{R}[x_C]_{r} \right\}, \\
  \Sigma_{\infty, C} &:=& \bigcup_{r\ge 0}\Sigma_{r, C}. 
  \end{eqnarray*}
  Note that if $C=\{1, \ldots, n\}$, then we have $\Sigma_{r, C}=\Sigma_r$ and $\Sigma_{\infty, C} = \Sigma$.
The sparse SDP relaxation problem with relaxation order $r$ for POP (\ref{eq:POP0}) is obtained from the following SOS relaxation problem: 
\begin{equation}\label{sparseSDPrelaxation}
\rho_r^{\mbox{\scriptsize sparse}} := \sup\left\{
\rho :
\begin{array}{l}
f - \rho = \sum_{h=1}^p\sigma_{0, h} + \sum_{j=1}^m \sigma_{j} f_j, \\
\sigma_{0, h}\in\Sigma_{r, C_h} \ (h=1, \ldots, p), \sigma_{j}\in\Sigma_{r_j, D_j} \ (j=1, \ldots, m)
\end{array}
\right\}, 
\end{equation}
where  $D_j$ is the union of some of the maximal cliques $C_1, \ldots, C_p$ such that $F_j\subseteq C_h$ and $r_j = r-\lceil\deg(f_j)/2\rceil$ for $j=1, \ldots, m$.

  It should be noted that {other} sparse SDP relaxations are proposed in \cite{Grimm07, Lasserre06b, Laurent09} and the asymptotic convergence  is proved. In contrast, the convergence of the sparse SDP relaxation (\ref{sparseSDPrelaxation}) is not shown in \cite{WAKI06}. 

We give an extension of Theorem \ref{th:main} to POP with correlative sparsity. If $C_1, \ldots, C_p\subseteq\{1, \ldots, n\}$ satisfy the following property, we refer this property as {\it the running intersection property} (RIP):
   \[
\forall h\in\{1, \ldots, p-1\}, \exists t\in\{1, \ldots, p\} \mbox{ such that  } C_{h+1}\cap (C_1\cup\cdots\cup C_h)\subsetneq C_t. 
 \]
For $C_1, \ldots, C_p\subseteq\{1, \ldots, n\}$, we define sets $J_1, \ldots, J_p$ as follows:
  \[
  J_h :=\left\{
  j\in\{1, \ldots, m\} :
  f_j\in\mathbb{R}[x_{C_h}]
  \right\}.
  \]
  Clearly, we have $\cup_{h=1}^p J_h = \{1, \ldots, m\}$.  In addition, we define
  \begin{eqnarray*}
  \psi_{r, h}(x) &:=& -\sum_{j\in J_h}f_j(x)\left(
  1 -\frac{f_j(x)}{R_j}
  \right)^{2r},\\
  \Theta_{r, h, b}(x) &:=&1 + \sum_{i\in C_h}\left(\frac{x_i}{b}\right)^{2r} 
  \end{eqnarray*}
  for $h=1, \ldots, p$. 

Using a proof similar to the one for the  theorem on convergence of the sparse SDP relaxation given in \cite{Grimm07}, we can establish the correlatively sparse case of Theorem \ref{th:main}. Indeed, we can obtain the theorem by using \cite[Lemma 4]{Grimm07} and Theorem \ref{th:main}. 

\begin{theorem}\label{th:sparsemain}
Assume that nonempty sets $C_1, \ldots, C_p\subseteq\{1, \ldots, n\}$ satisfy (RIP) and we can decompose $f$ into $f= \hat{f}_1+\cdots +\hat{f}_p$ with $\hat{f}_h\in\mathbb{R}[x_{C_h}] \ (h=1, \ldots, p)$. Under the assumptions of Theorem \ref{th:main}, there exists $\tilde{r}\in\mathbb{N}$ such that for all $r\ge\tilde{r}$, $f - \rho + \sum_{h=1}^p\psi_{r, h}$ is positive over $B=[-b, b]^n$. In addition, for every $\epsilon>0$, there exists $\hat{r}\in\mathbb{N}$ such that for all $r\ge \hat{r}$, 
\begin{equation} \label{eq:sparsemain}
f - \rho +\epsilon \sum_{h=1}^p\Theta_{r, h, b} + \sum_{h=1}^p\psi_{\tilde{r}, h}\in \Sigma_{\infty, C_1} +\cdots +\Sigma_{\infty, C_p}. 
\end{equation}
\end{theorem}
Note that if $p=1$, {\it i.e.}, $C_1 = \{1,\ldots,n\}$,
then we have $\psi_{r, 1} = \psi_{r}$ and $\Theta_{r, 1, b} = \Theta_{r, b}$, and thus Theorem \ref{th:sparsemain} is reduced to Theorem \ref{th:main}.
Therefore, we will concentrate our effort to prove Theorem \ref{th:sparsemain}
in the following. In addition, we remark that it would follow from \cite[Theorem 5]{Grimm07} that 
\eqref{eq:sparsemain} holds without the polynomial $\epsilon \sum_{h=1}^p\Theta_{r, h, b}$ if we assume that all  quadratic modules generated by $f_j \ (j\in C_h)$ for all $h=1, \ldots, p$ are \redch{A}{a}rchimedean. 

To prove Theorem \ref{th:sparsemain}, we use Lemma 4 in \cite{Grimm07}
and Corollary 3.3 of \cite{LaNe07}.

\begin{lemma}\label{grimmlemma}(modified version of \cite[Lemma 4]{Grimm07}) 
Assume that we decompose $f$ into $f=\hat{f}_1+\cdots +\hat{f}_p$ with $\hat{f}_h\in\mathbb{R}[x_{C_h}]$ and $f > 0$ on $K$. Then, for any bounded set $B\subseteq\mathbb{R}^n$, there exist 
$\tilde{r}\in\mathbb{N}$ and $g_h\in\mathbb{R}[x_{C_h}]$ with $g_h > 0$ on $B$ such that for every $r\ge\tilde{r}$, 
\[
f = -\sum_{h=1}^p\psi_{r, h} + \sum_{h=1}^p g_h.
\]
\end{lemma}
\begin{remark}\label{Lemma}
The original statement in \cite[Lemma 4]{Grimm07} is slightly different from Lemma \ref{grimmlemma}. In \cite[Lemma 4]{Grimm07}, it is proved  that there exists $\lambda\in (0, 1]$, $\tilde{r}\in\mathbb{N}$ and $g_h\in\mathbb{R}[x_{C_h}]$ with $g_h > 0$ on $B$ such that 
\[
f = \sum_{h=1}^p\sum_{j\in J_h}\left(1 - \lambda f_j\right)^{2\tilde{r}}f_j + \sum_{h=1}^p g_h.
\]
In Appendix \ref{sec:appendix}, we {establish} the correctness of Lemma \ref{grimmlemma} by using \cite[Lemma 4]{Grimm07}. 
\end{remark} 

\begin{lemma}{(Corollary 3.3 of \cite{LaNe07})} \label{lem:3}
Let $f\in \R[x]$ be a polynomial nonnegative on $[-1,1]^n$.
For arbitrary $\epsilon>0$, there exists some $\hat{r}$ such that
for every $r\geq \hat{r}$, the polynomial $(f+\epsilon \Theta_r)$ is a 
SOS.
\end{lemma}
\begin{pf*}{Theorem \ref{th:sparsemain}}
We may choose $[-b, b]^n$ as $B$ in Lemma \ref{grimmlemma}. It follows from the assumption in Theorem  \ref{th:sparsemain} that we can decompose $f-\rho$ into $(\hat{f}_1-\rho) + \hat{f}_2 + \cdots + \hat{f}_p$.  Since $\hat{f}_1-\rho\in\mathbb{R}[x_{C_1}]$,  it follows from Lemma \ref{grimmlemma} that there exists $\tilde{r}\in\mathbb{N}$ and $g_h\in\mathbb{R}[x_{C_h}]$ with $g_h > 0$ on $B$ such that for every $r\ge \tilde{r}$, 
\[
f -\rho = (\hat{f}_1-\rho) +\hat{f}_2+ \cdots + \hat{f}_p = -\sum_{h=1}^p\psi_{r, h} + \sum_{h=1}^p g_h.
\]
Therefore, the polynomial $f-\rho + \sum_{h=1}^p\psi_{r, h}$ is positive on $B$ for all $r\ge \tilde{r}$. 

For simplicity, we fix $h$ and  define $C_h=\{c_1, \ldots, c_k\}$. Then, $g_h$ consists of the $k$ variables $x_{c_1}, \ldots, x_{c_k}$. 
Since $g_h >0$ on $B$, it is also positive on $B^{\prime}:=\{(x_{c_1}, \ldots, x_{c_k}): -b\le x_{c_j} \le b \ (j=1, \ldots, k)\}$.  We define $\hat{g}_h(y)= g_h(by)$.  Since $g_h$ is positive on $B^{\prime}$, $\hat{g}_h\in\mathbb{R}[y_{c_1}, \ldots, y_{c_k}]$ is also positive on the set $\{(y_{c_1}, \ldots, y_{c_k}):  -1\le y_{c_j} \le 1 \ (j=1, \ldots, k)\}$.  Applying Lemma \ref{lem:3} to $\hat{g}_h$, for all $\epsilon>0$, there exists  $\hat{r}_h\in\mathbb{N}$ such that for every $r\ge \hat{r}_h$, 
\[
 \hat{g}_h(y_{c_1}, \ldots, y_{c_k}) + \epsilon\sum_{i=1}^k y_{c_i}^{2r} =\sigma_h(y_{c_1}, \ldots, y_{c_k})
\]
for some $\sigma_h\in \Sigma_{\infty, C_h}$. Substituting $x_{c_1}=by_{c_1}, \ldots, x_{c_k}=by_{c_k}$, we obtain
\[
g_h +\epsilon\Theta_{r, h, b} \in\Sigma_{\infty, C_h}. 
\]

We fix $\epsilon>0$. Applying the above discussion to all $h=1, \ldots, p$, we obtain the numbers $\hat{r}_1, \ldots, \hat{r}_p$. We denote the maximum over $\hat{r}_1, \ldots, \hat{r}_p$ by $\hat{r}$. Then, we have  
\[
f-\rho +\epsilon\sum_{h=1}^p\Theta_{r, h, b} + \sum_{h=1}^p\psi_{\tilde{r}, h}\in \Sigma_{\infty, C_1} +\cdots +\Sigma_{\infty, C_p}
\]
for every $r\ge \hat{r}$. 
\end{pf*}


\subsection{Extension to POP with symmetric cones}\label{subsec:symcone}

In this subsection, we extend Theorem \ref{th:main} to POP over  symmetric cones, {\it i.e.}, 
\begin{equation}\label{POPsymcone}
f^*:=\inf_{x\in\mathbb{R}^n}\left\{
f(x) : G(x)\in\mathcal{E}_+
\right\}, 
\end{equation}
where $f\in\mathbb{R}[x]$, $\mathcal{E}_+$ is a symmetric cone associated with an $N$-dimensional Euclidean Jordan algebra $\mathcal{E}$, and $G$ is $\mathcal{E}$-valued polynomial in $x$. The feasible region $K$ of POP (\ref{POPsymcone}) is $\{x\in\mathbb{R}^n : G(x)\in\mathcal{E}_+\}$. Note that if $\mathcal{E}$ is $\mathbb{R}^m$ and $\mathcal{E}_+$ is the nonnegative orthant $\mathbb{R}^m_+$, then (\ref{POPsymcone}) is identical to (\ref{eq:POP0}). In addition, $\mathbb{S}^n_+$, the cone of $n\times n$ symmetric positive semidefinite matrices,  is a symmetric cone, the bilinear matrix inequalities can be formulated as (\ref{POPsymcone}). 

To construct $\psi_r$ for (\ref{POPsymcone}), we introduce some notation and symbols. The Jordan product and inner product of $x, y\in\mathcal{E}$ are denoted by, respectively,  $x\circ y$ and $x\bullet y$. Let $e$ be the identity element in the Jordan algebra $\mathcal{E}$. For any $x\in\mathcal{E}$, we have $e\circ x = x\circ e = x$. We can define eigenvalues for all elements in the Jordan algebra $\mathcal{E}$, generalizing those for Hermitian matrices. See \cite{Faruat} for the details.  We construct $\psi_r$ for (\ref{POPsymcone}) as follows: 
\begin{eqnarray}
\nonumber
M&:=& \sup\left\{
\mbox{maximum absolute eigenvalue of }G(x):  x\in \bar{K}
\right\}, \\
\label{psi_for_symcone}
\psi_r(x) &:=& -G(x)\bullet \left(
e - \frac{G(x)}{M}
\right)^{2r}, 
\end{eqnarray}
where we define $x^k:= x^{k-1}\circ x$ for $k\in\mathbb{N}$ and $x\in\mathcal{E}$. 

Lemma 4 in \cite{Kojima07} shows that $\psi_r$ defined in (\ref{psi_for_symcone}) has the same properties as $\psi_r$ in Theorem \ref{th:main}. 

\begin{theorem}\label{th:POPsymcone}
For a given $\rho$, suppose that $f(x)-\rho >0$ for every $x\in\bar{K}$. Then, there exists $\tilde{r}\in\mathbb{N}$ such that for all $r\ge\tilde{r}$, $f -\rho +\psi_r$ is positive over $B$.  Moreover, for any $\epsilon >0$, there exists $\hat{r}\in\mathbb{N}$ such that for every $r\ge\hat{r}$, 
\[
f -\rho + \epsilon\Theta_{r, b} + \psi_{\tilde{r}} \in\Sigma. 
\]

\end{theorem}

\subsection{Another perturbed sums of squares theorem}\label{subsec:alternative}
In this subsection, we present another perturbed sums of squares theorem for POP (\ref{eq:POP0}) which is obtained by combining results in \cite{Kim05, Lasserre06}. 

To use the result in \cite{Kim05}, we introduce some notation and symbols. We assume that $K\subseteq B:=[-b, b]^n$. We choose $\gamma\ge 1$ such that for all $j=0, 1, \ldots, m$, 
\begin{eqnarray*}
|f_j(x)/\gamma| &\le& 1 \mbox{ if } \|x\|_{\infty}\le\sqrt{2}b, \\
|f_j(x)/\gamma| &\le& \|x/b\|_{\infty}^d \mbox{ if } \|x\|_{\infty}\ge\sqrt{2}b, \\
\end{eqnarray*}
where $f_0$ denotes the objective function $f$ in POP (\ref{eq:POP0}),  and $d =\max\{\deg(f), \deg(f_1), \ldots, \deg(f_m)\}$. For $r\in\mathbb{N}$, we define 
\begin{eqnarray*}
\psi_r(x) &:=& -\sum_{j=1}^m\left(1-\frac{f_j(x)}{\gamma}\right)^{2r} f_j(x), \\
 \phi_{r, b}(x) &:=&-\frac{(m+2)\gamma}{b^2}\sum_{i=1}^n\left(
 \frac{x_i}{b}
 \right)^{2d(r+1)}(b^2 -x_i^2). 
\end{eqnarray*}

From (a), (b) and (c) of Lemma 3.2 in \cite{Kim05}, we obtain the following result:
\begin{proposition}\label{GLFpos} 
Assume that the feasible region $K$ of POP (\ref{eq:POP0}) is contained in $B=[-b, b]^n$. In addition, we assume that for $\rho\in\mathbb{R}$, we have $f-\rho > 0$ over $K$. 
Then there exists $\tilde{r}\in\mathbb{N}$ such that for all $r\ge\tilde{r}$, $(f-\rho +\psi_r +\phi_{r, b})$ is positive over $\mathbb{R}^n$.
\end{proposition}

We remark that we do not need to impose the assumption on the compactness of $K$ in Proposition \ref{GLFpos}. Indeed, we can drop it by replacing $K$ by $\bar{K}$ defined in Subsection \ref{subsec:main}  as in Theorem \ref{th:main}.

Next, we describe a result from \cite{Lasserre06} which is useful in deriving  another perturbed sums of squares theorem. 
\begin{theorem} ((iii) of Theorem 4.1 in \cite{Lasserre06})\label{lasserreTh}
Let $f\in\mathbb{R}[x]$ be a nonnegative polynomial. Then for every $\epsilon>0$, there exists $\hat{r}\in\mathbb{N}$ such that for all $r\ge\hat{r}$, 
\[
f + \epsilon\theta_r\in\Sigma, 
\]
where $\theta_r(x) := \sum_{i=1}^n\sum_{k=0}^r(x_i^{2k}/k!)$. 
\end{theorem}

By incorporating Proposition \ref{GLFpos} with Theorem \ref{lasserreTh}, we obtain yet another perturbation theorem. 

\begin{theorem}\label{anotherPT}
We assume that for $\rho\in\mathbb{R}$, we have $f-\rho > 0$ over $K$. Then we have 
\begin{enumerate}
\item there exists $\tilde{r}\in\mathbb{N}$ such that for all $r\ge\tilde{r}$, $(f-\rho +\psi_r +\phi_{r, b})$ is positive over $\mathbb{R}^n$; 
\item moreover, for every $\epsilon>0$, there exists $\hat{r}\in\mathbb{N}$ such that for all $r\ge\hat{r}$, 
\[
(f -\rho +\psi_{\tilde{r}} +\phi_{\tilde{r}, b}+\epsilon\theta_r) \in\Sigma. 
\]
\end{enumerate}
\end{theorem}

We give an SDP relaxation  analogous to (\ref{eq:sparse1}), based on Theorem \ref{anotherPT}, as follows:
\begin{equation}\label{anotherSDPr}
\eta(\epsilon, \tilde{r}, r) :=\sup\left\{
\eta :
\begin{array}{l}
f -\eta +\epsilon\theta_r -\displaystyle\sum_{j=1}^m f_j\sigma_j -\displaystyle\sum_{i=1}^n (b^2 -x_i^2)\mu_i = \sigma_0, \\
\sigma_0\in\Sigma_r, \sigma_j\in\Sigma(\tilde{r}\tilde{\mathcal{F}}_j), \mu_i\in\Sigma(\{d(\tilde{r}+1)e_i\})
\end{array}
\right\}, 
\end{equation}
for some $r\ge\tilde{r}$, where $e_i$ is the $i$th standard unit vector in $\mathbb{R}^n$. One of the differences between (\ref{eq:sparse1}) and (\ref{anotherSDPr}) is that (\ref{anotherSDPr}) has $n$ SOS variables $\mu_1, \ldots, \mu_n$. These variables correspond to nonnegative variables in the SDP formulation, but not positive semidefinite matrices,  since these consist of a single monomial. On the other hand, it is difficult to estimate $\tilde{r}$ in the SDP relaxations (\ref{eq:sparse1}) and (\ref{anotherSDPr}), and thus we could not compare the size and the quality of the optimal value of (\ref{eq:sparse1}) with (\ref{anotherSDPr}) so far.


We obtain a result similar technique to Theorem \ref{th:conv}. We omit the proof because we obtain the inequalities by applying a proof similar to that of Theorem \ref{th:conv}.
\begin{theorem}
For every $\epsilon>0$, there exists $r, \tilde{r}\in\mathbb{N}$ such that 
$f^*-\epsilon\le \eta(\epsilon, \tilde{r}, r)\le f^* + \epsilon n e^{b^2}$. 
\end{theorem}

\section{Concluding Remarks} \label{sec:conclusion}

We mention other {research} related to our work related to Theorem \ref{th:main}.  A common element in all of these approaches is to use perturbations $\epsilon\theta_r(x)$ or $\epsilon\Theta_r(x)$ for finding an approximate solution of a given POP. 

In \cite{HANZON03, JibeteanPHD},   the authors added $\epsilon\Theta_r(x)$ to the objective function of a given unconstrained POP and used algebraic techniques to find a solution. In \cite{Jibetean05},  the following equality constraints were added in the perturbed unconstrained POP and Lasserre's SDP relaxation was applied to the new POP:  
\[
 \frac{\partial f_0}{\partial x_i} +  2r\epsilon x_i^{2r-1} = 0 \ (i=1, \ldots, n). 
\]

Lasserre in \cite{Lasserre07} proposed an SDP relaxation via $\theta_r(x)$ defined in Theorem \ref{lasserreTh} and a perturbation theorem for semi-algebraic set defined by equality constraints $g_k(x) = 0$ \ $(k=1, \ldots, m)$. The SDP relaxation can be applied to the following equality constrained POP:
\begin{equation}\label{POP2}
\inf_{x\in\mathbb{R}^n}\left\{
f_0(x) : 
g_k(x) = 0 \ (k=1, \ldots, m)
\right\}; 
\end{equation}
To obtain the SDP relaxations, $\epsilon\theta_r(x)$ is added
to the objective function in POP (\ref{POP2})  and 
the equality constraints in POP (\ref{POP2}) is replaced 
by $g_k^2(x)\le 0$.
In the resulting SDP relaxations, $\theta_r(x)$ is explicitly introduced and variables associated with constraints $g_k^2(x)\le 0$ are not positive semidefinite matrices, but nonnegative variables.  

\medskip

In this paper, we {present} a perturbed SOS theorem (Theorem \ref{th:main})
and its extensions, and propose a new sparse relaxation called 
Adaptive SOS relaxation. During the course {of the paper}, we have {shed some light on} why 
Lasserre's SDP relaxation calculates the optimal value of POP
even if its SDP relaxation has a different optimal value.
The numerical experiments clearly show that Adaptive SOS relaxation 
is promising, justifying the need for future research in this direction.

Of course, if the original POP is dense, i.e., %
$\tilde{F}_j$ contains many elements for almost all $j$, 
then the proposed relaxation has little effect in reducing 
the SDP relaxation. {However,} in real {applications}, such {cases seem} rare.

In the numerical experiments, we sometimes observe that 
the behaviors of SeDuMi and SDPT3 are very different each other.
See, for example, Table \ref{subsec5.3_table}. In the column of Adaptive SOS,
SeDuMi solved significantly {fewer} problems than SDPT3.
On the other hand, there are several cases where SeDuMi 
outperforms SDPT3.
For such an example, see the sparse relaxation column of Table \ref{subsec5.4_table_density0.6}.
This is why we present the results of both solvers in every table.
In solving a real problem, {one} should be {very} careful in choosing 
the appropriate SDP solver for the problem {at hand}.

\section*{Acknowledgements}
The first author was supported by in part by a Grant-in-Aid for Scientific Research (C) 19560063. 
The second author was supported in part bya Grant-in-Aid for Young Scientists (B) 22740056. 
The third author was supported in part by a Discovery Grant from NSERC, a research grant from University of Waterloo and by ONR research grant N00014-12-10049. 

\appendix
\section{A proof of Lemma \ref{grimmlemma}}\label{sec:appendix}

As we have already mentioned in Remark \ref{Lemma}, Lemma \ref{grimmlemma} is slightly different from the original one in \cite[Lemma 4]{Grimm07}. To show the correctness of Lemma \ref{grimmlemma}, we use the following lemma:

\begin{lemma}(\cite[Lemma 3]{Grimm07})\label{grimmlemma3}
Let $B\subseteq\mathbb{R}^n$ be a compact set. Assume that nonempty sets $C_1, \ldots, C_p\subseteq\{1, \ldots, n\}$ satisfy (RIP) and we can decompose $f$ into $f= \hat{f}_1+\cdots +\hat{f}_p$ with $\hat{f}_h\in\mathbb{R}[x_{C_h}] \ (h=1, \ldots, p)$. In addition,  suppose that $f >0$ on $B$. Then there exists $g_h\in\mathbb{R}[x_{C_h}]$ with $g_h>0$ on $B$ such that
\[
f = g_1 +\cdots + g_p.
\] 
\end{lemma}

We can prove Lemma \ref{grimmlemma} in a manner similar to \cite[Lemma 4]{Grimm07}. We define $F_r: \mathbb{R}^n\to\mathbb{R}$ as follows:
\[
F_r = f - \sum_{h=1}^p\psi_{r, h}. 
\]
We recall that $\psi_{r, h} = \sum_{j\in C_h}(1-f_j/R_j)^{2r}f_j$ for all $h=1, \ldots, p$ and $r\in\mathbb{N}$, and that $R_j$ is the maximum {value} of $|f_j|$ on $B$ for all $j=1, \ldots, m$. It follows from the definitions of  $\psi_{r, h}$ and $R_j$ that we have $\psi_{r, h}\ge \psi_{r+1, h}$ on $B$  for all $h=1, \ldots, p$ and $r\in\mathbb{N}$, and thus we have $F_r\le F_{r+1}$ on $B$.  In addition, we can prove that (i) on $B\cap K$, $F_r\to f$ as $r\to\infty$, and (ii) on $B\setminus K$, $F_r\to\infty$ as $r\to\infty$. Since $B$ is compact, it follows from (i), (ii) and the positiveness of $f$ on $B$ that there exists $\tilde{r}\in\mathbb{N}$
 such that for every $r\ge \tilde{r}$, $F_{r} > 0$ on $B$. Applying Lemma \ref{grimmlemma3} to $F_r$, we obtain the desired result. 


\end{document}